\documentclass[12pt]{iopart}
\usepackage{epsfig,amsfonts}
\newcommand{\norm}[1]{\parallel\!\!#1\!\!\parallel}

\begin{document}
\title{ Coarse-grained dynamics of an activity bump 
	in a neural field model}
\author{C R Laing$^1$, T A Frewen$^2$ and I G Kevrekidis$^2$}
\address{$^1$ Institute of Information and Mathematical Sciences, Massey University,
Private Bag 102-904 NSMC, Auckland, New Zealand}
\address{$^2$ Department of Chemical Engineering,  Princeton University, NJ 08544-5263, USA}
\ead{c.r.laing@massey.ac.nz}

\begin{abstract}

We study a stochastic nonlocal PDE, arising in the context of modelling spatially distributed neural
activity, which is capable of sustaining stationary and moving spatially-localized
``activity bumps''.
This system is known
to undergo a
pitchfork bifurcation in bump speed as a parameter (the strength of
adaptation) is changed; yet increasing the noise intensity
effectively slowed the motion of the bump.
Here we revisit the system from the point of view of
describing the high-dimensional stochastic dynamics in terms of the
effective dynamics of a single scalar ``coarse" variable.
We show that such a reduced description in the form of an effective Langevin
equation characterized by a double-well potential
is quantitatively successful.
The effective potential can be extracted using short, appropriately-initialized bursts of
direct simulation.
We demonstrate this approach in terms of (a) an experience-based ``intelligent" choice
of the coarse observable and (b) an observable obtained through data-mining direct
simulation results, using a
diffusion map approach.

\end{abstract}
\ams{60H15,65P30,37G35,45K05,82C31}
\maketitle

\section{Introduction}
\label{sect:intro}

Pattern formation at large spatial scales in the cortex is a field of much recent
interest~\cite{benbar95,cooowe04,erm98,pinerm01A}. Neural dynamics are intrinsically noisy,
so it is natural to study the effects of noise on models for such pattern formation.
In this paper we revisit a noisy model
for pattern formation in a neural system~\cite{lailon01} but this time concentrating on describing
its dynamics in a low-dimensional way.
We will study the partial integrodifferential system
\begin{eqnarray}
   \frac{\partial u(x,t)}{\partial t} = -u(x,t)+\int_{-\pi}^{\pi}J(x-y)f[I+u(y,t)-a(x,t)]dy \label{eq:dudt} \\
   \tau\frac{\partial a(x,t)}{\partial t} = Au(x,t)-a(x,t) \label{eq:dadt}
\end{eqnarray}
on the domain $-\pi\leq x\leq\pi$ with periodic boundary conditions,
originally analyzed by Laing and Longtin~\cite{lailon01} as a simple model for pattern formation
in a neural field.  The variable $u(x,t)$ describes the neural activity at position $x$ and time $t$.
The function $J$ describes the spatial coupling within the network, with $J(x-y)$ being the
strength of coupling between neurons at position $x$ and neurons at position $y$.
The function $f(I)$
describes the firing rate of a single neuron with input $I$.
The quantities $A,I$ and $\tau$ are constant. The value of $A$ determines the strength of adaptation,
$I$ is a background current applied to all neurons and $\tau$ is the adaptation time constant.
Equation~(\ref{eq:dadt}) represents the effects of spike frequency adaptation.
Without~(\ref{eq:dadt}) systems of this form have been used to model
orientation tuning in the visual system,
working memory and the head direction  system~\cite{erm98,benbar95}.
Variants including~(\ref{eq:dadt}) have been subsequently explored
~\cite{cooowe04,cooowe05,pinerm01A}.

In our numerical implementation
we spatially discretise the domain with $M$ equally spaced points, resulting in the $2M$ ODEs
\begin{eqnarray}
   \frac{du_i}{dt} = -u_i+\frac{2\pi}{M}\sum_{j=1}^{M}J_{ij}f(I+u_i-a_i) \label{eq:dudtA}\\
   \tau\frac{da_i}{dt} = Au_i-a_i \label{eq:dadtA}
\end{eqnarray}
for $i=1,\ldots, M$, where $J_{ij}=J(2\pi|i-j|/M)$.
As in Laing and Longtin~\cite{lailon01}
we set $I=-0.1,\tau=5$ and use $f(x)=(1+\mbox{tanh}(10x))/2$
and $J(x)=0.05+0.24\cos{x}$, with $M=100$.
We add noise to the system by adding a
Gaussian white noise term
$\xi_i(t)$ to each of~(\ref{eq:dudtA}), with $\langle\xi_i(t)\rangle=0$ and
$\langle\xi_i(t)\xi_j(s)\rangle=2\eta\nu_{ij}\delta(t-s)$, where $\nu_{ij}=0$ if $i\neq j$ and
1 if $i=j$.
Initially we set the parameters $A=0.17$ and $\eta=10^{-4}$; we will study the effect of
varying them below.

\begin{figure}
\leavevmode
\epsfxsize=5.9in
\epsfbox{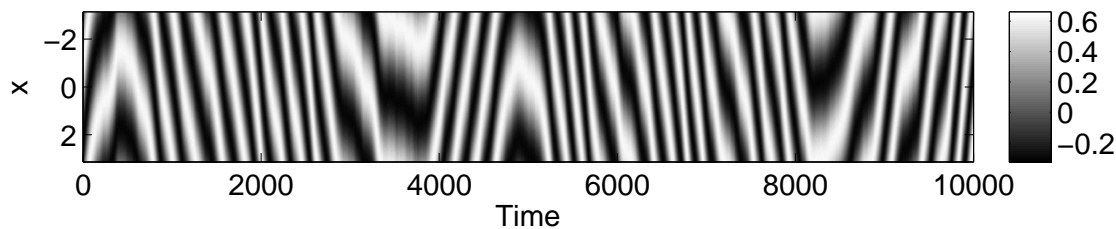}
\caption{A typical simulation of the system~(\ref{eq:dudtA})-(\ref{eq:dadtA}).
The profile of $u$ is shown as a function of space and time (values indicated by shadebar). Parameters are
$M=100,\eta=10^{-4},A=0.17,\tau=5.$}
\label{fig:example}
\end{figure}

At the deterministic limit $\eta=0$, for these parameter values the system dynamics ultimately exhibit
a single stable activity bump that moves around the (periodic) domain in
one direction (determined by the initial condition), at a constant velocity.
If $\eta\neq 0$,
the bump motion becomes stochastic, exhibiting occasional switches in direction.
This effect, studied by Laing and Longtin~\cite{lailon01}, is illustrated in figure~\ref{fig:example}; they concentrated on the net distance
travelled during a finite time interval.
For small $\eta$ the bump moved in only one direction,
so that the net distance travelled was quite far (effectively proportional
to the time);  for $\eta$ greater than a certain value, switches
in direction lead to a much smaller net distance travelled.
This phenomenon was referred  to as
{\em noise-induced stabilisation}.

Laing and Longtin~\cite{lailon01} showed that for a coupling function of the form used here,
the
system~(\ref{eq:dudt})-(\ref{eq:dadt}) underwent a pitchfork bifurcation in speed as $A$
was increased.
They then added Gaussian white noise to the normal form of a pitchfork
bifurcation, and qualitatively reproduced the noise-induced stabilisation.
However, the actual derivation of a ``noisy normal form" from the original system with noise
was not attempted.
In this paper we will demonstrate numerically that the noisy discretised
system~(\ref{eq:dudtA})-(\ref{eq:dadtA}) undergoes what we will characterize as
an {\em effective pitchfork bifurcation} in speed as $A$
is increased, and investigate the effects of varying noise intensity.
This computer-assisted
bifurcation analysis will be performed assuming that the dynamics of the high-dimensional
stochastic system~(\ref{eq:dudtA})-(\ref{eq:dadtA}) can be effectively described by the dynamics
of a single, scalar coarse-grained {\em observable} (variable).
The evolution of this variable can be determined by running many short, appropriately-initialized
simulations of~(\ref{eq:dudtA})-(\ref{eq:dadtA}), in the spirit of the
``equation-free'' framework~\cite{kevgea03,lai06}.

The value of our scalar variable $V$ (defined below) can be interpreted as the position of a
``particle'' moving in a potential, subject to noise.
Beyond the effective pitchfork bifurcation
the potential is double-welled; before it, it is single-welled.
Each of the two wells corresponds to persistent motion in one direction (left or right);
in the case of a single well we have a stationary (on average) activity bump.
Since the system is isotropic, we expect the double well to be symmetric.
Although we cannot analytically derive the effective potential and noise in our
hypothesized Langevin description for $V$, we will show how they can be estimated using
appropriately initialized short bursts of simulation of~(\ref{eq:dudtA})-(\ref{eq:dadtA}).
Once the effective potential is approximated, we can locate its maxima and minima
(which take the place of unstable and stable fixed points for the particle, respectively)
and estimate average transition times between minima (corresponding to bump direction
switches).
By determining how the potential changes as $A$ and $\eta$ are varied we can also numerically
obtain a bifurcation diagram for its extrema, which takes the place of the traditional
steady-state bifurcation diagram for the deterministic problem; identifying transition states
(potential maxima, corresponding to unstable steady states) would be extremely difficult
if not impossible using long simulations alone.
An interesting point is that we can quantify the relationship between the noise intensity in the original
system, $\eta$, and the effective noise intensity, $D(V)$, to which the particle in the potential is
subject.

\section{Results}

\subsection{Reconstructing the potential}
\label{subsect:reconstr}

\begin{figure}
\leavevmode
\epsfxsize=5.9in
\epsfbox{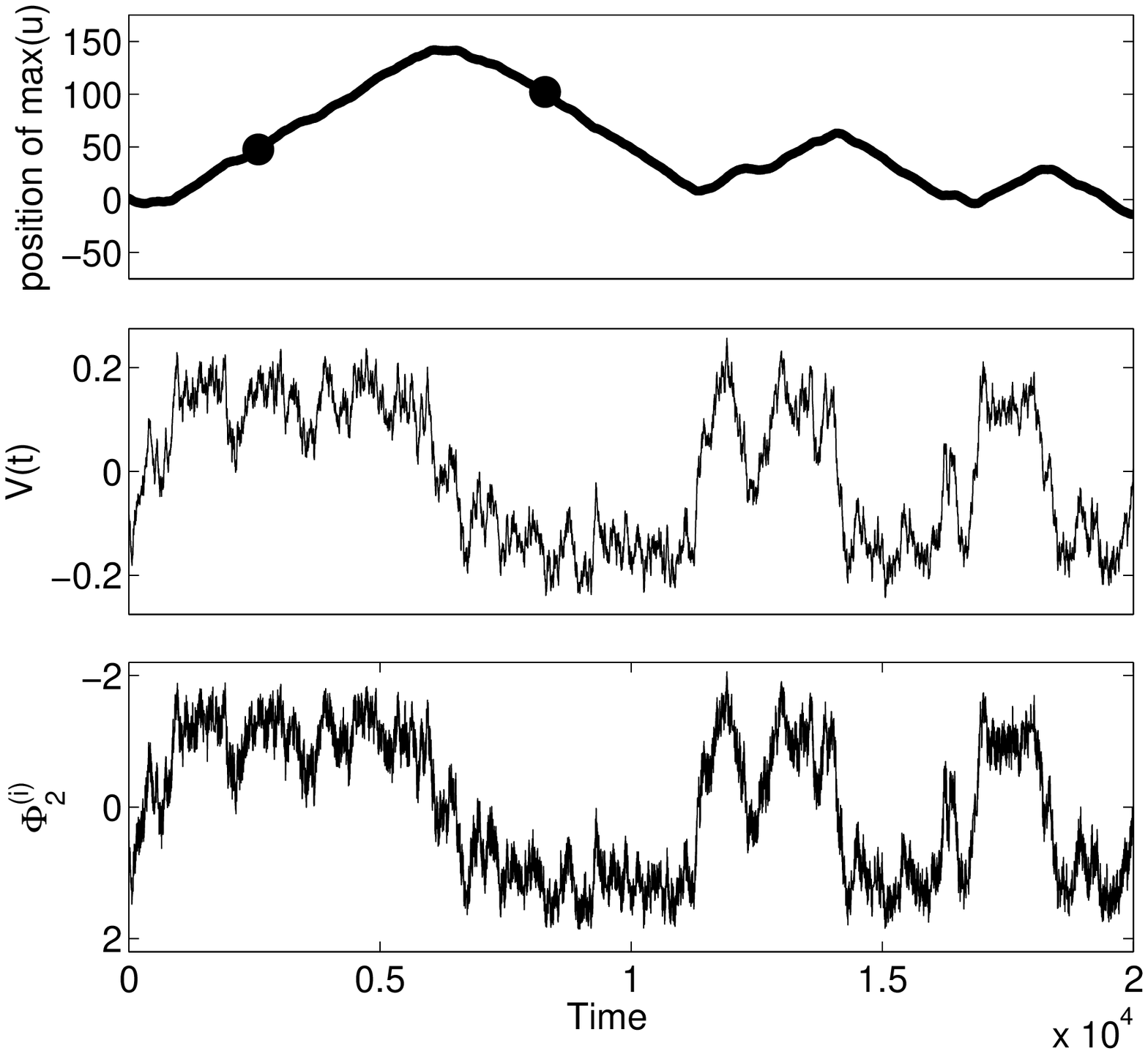}
\epsfxsize=5.in
\centerline{
\epsfbox{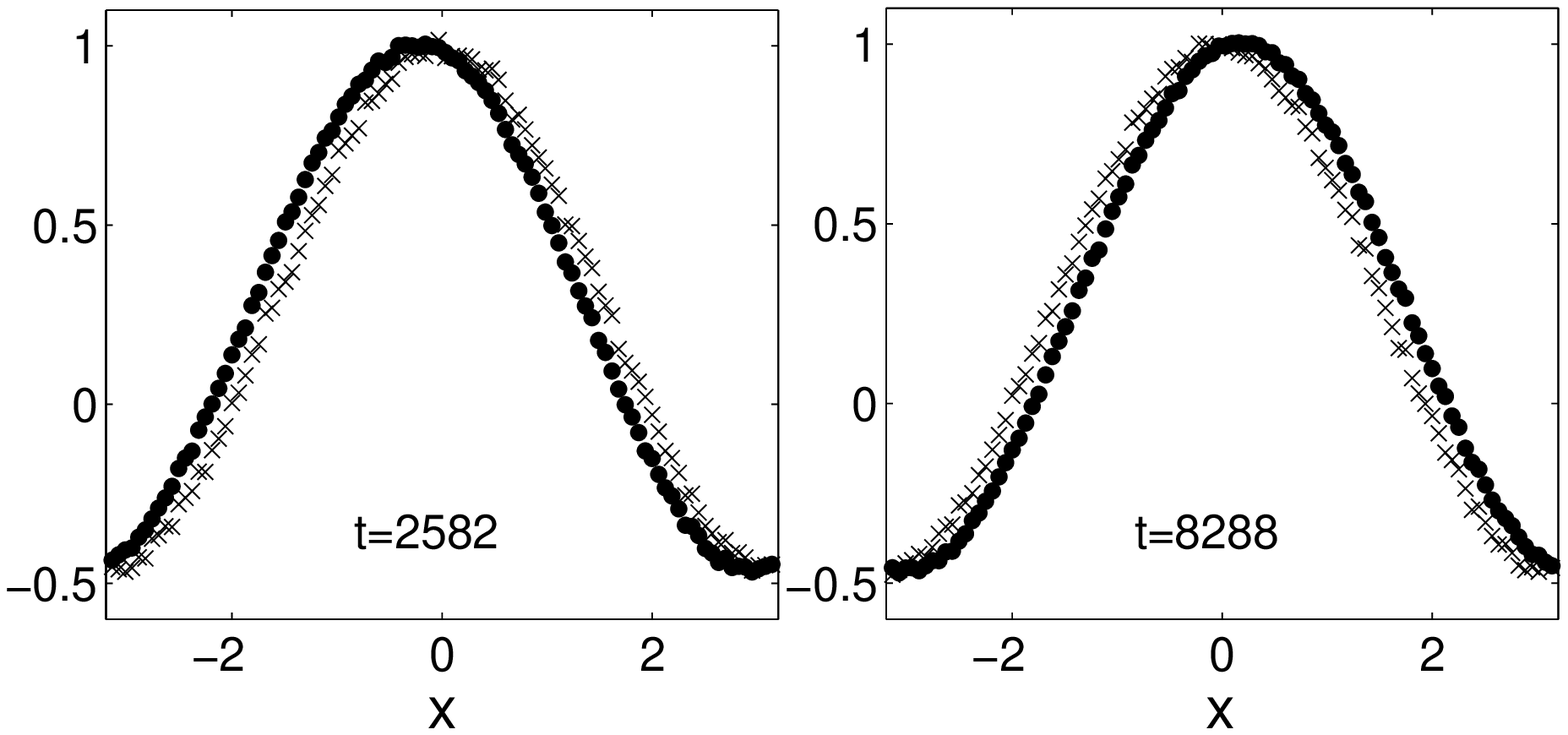}
}
\caption{Simulation results for system~(\ref{eq:dudtA})-(\ref{eq:dadtA}).
Top row: position of the peak in the $u$ profile as a function of time; second
and third rows: $V(t)$  and $\Phi_2^{(i)}$ (defined in the text),
respectively for the same time interval; bottom row: spatial profiles
of $u$ (cross) and $a$ (filled circles) at times marked in top panel. Parameters are the same as in
figure~\ref{fig:example}.}
\label{fig:resmulti}
\end{figure}

\begin{figure}
\leavevmode
\epsfxsize=5.9in
\epsfbox{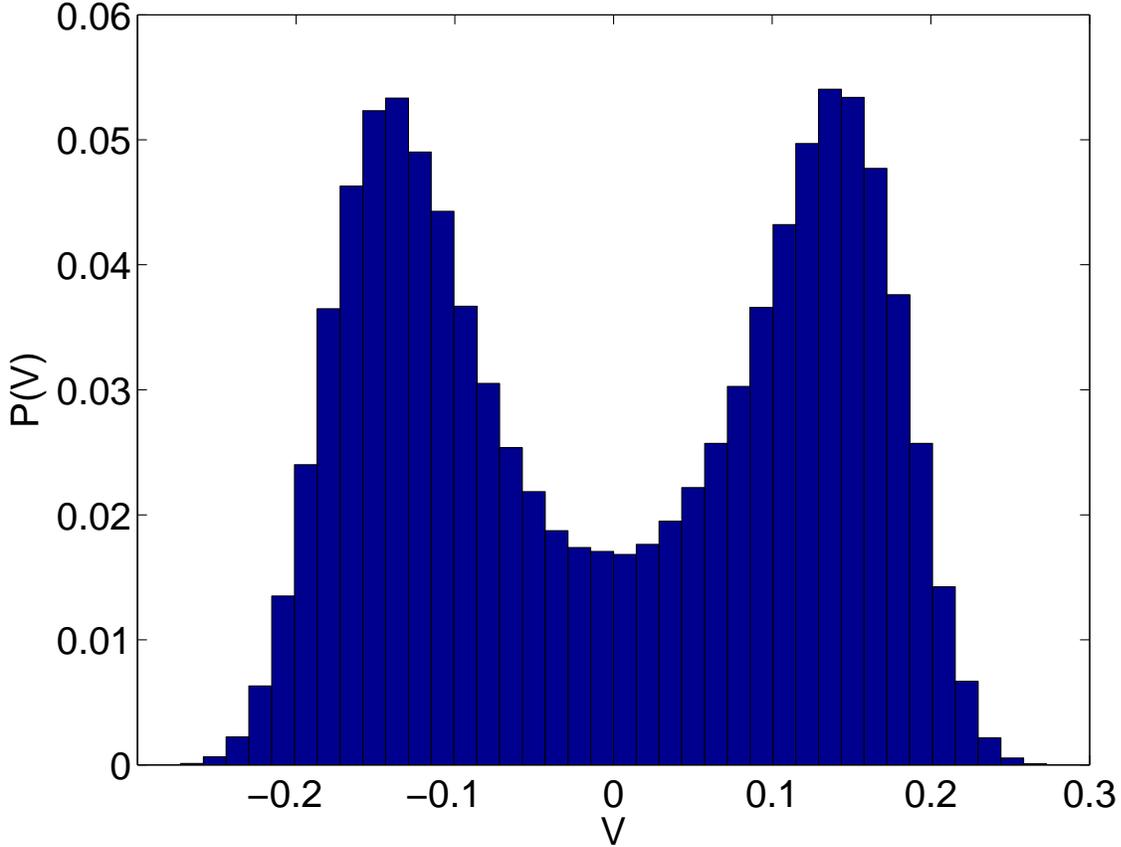}
\caption{A typical distribution of values of $V$. Parameter values are $A=0.17$ and
$\eta=10^{-4}$.}
\label{fig:hist}
\end{figure}

We choose as our coarse scalar variable, $V(t)$,
the instantaneous difference in position between the peak of $u(x,t)$ and the peak of $a(x,t)$.
This choice comes from observations of the system dynamics: practically all system profiles
during a long simulation are characterized by a single peak in each of the variable fields.
Furthermore, the difference between these peaks is clearly related to the instantaneous
speed of the bump:  the peak of $a$ normally lags the peak of $u$ relative to the direction
of bump motion.
One can clearly see in figure~\ref{fig:resmulti} that periods of bump motion to the left (right) are characterized
by effectively constant negative (positive) $V(t)$.
Factoring out the small amplitude noisy oscillations riding on the bump waveform,
we identify the position of the peak of $u(x,t)$ as the location $c_u$ that satisfies
\begin{equation}
   \sum_{i=1}^M\sin{(x_i-c_u)}u_i(t)=0
\end{equation}
where $x_i=2\pi i/M$.
(Effectively, we are finding the amount by which a simple sine wave must be shifted
so that its inner product with $u(x,t)$ is zero.)
Similarly, $c_a$ satisfies
\begin{equation}
   \sum_{i=1}^M\sin{(x_i-c_a)}a_i(t)=0
\end{equation}
so that  $V\equiv c_u-c_a$.
Results for a typical simulation are shown in
figure~\ref{fig:resmulti} and
a typical probability distribution of $V$, for the same parameters, is
shown in figure~\ref{fig:hist}.

\begin{figure}
\leavevmode
\epsfxsize=5.0in
\epsfbox{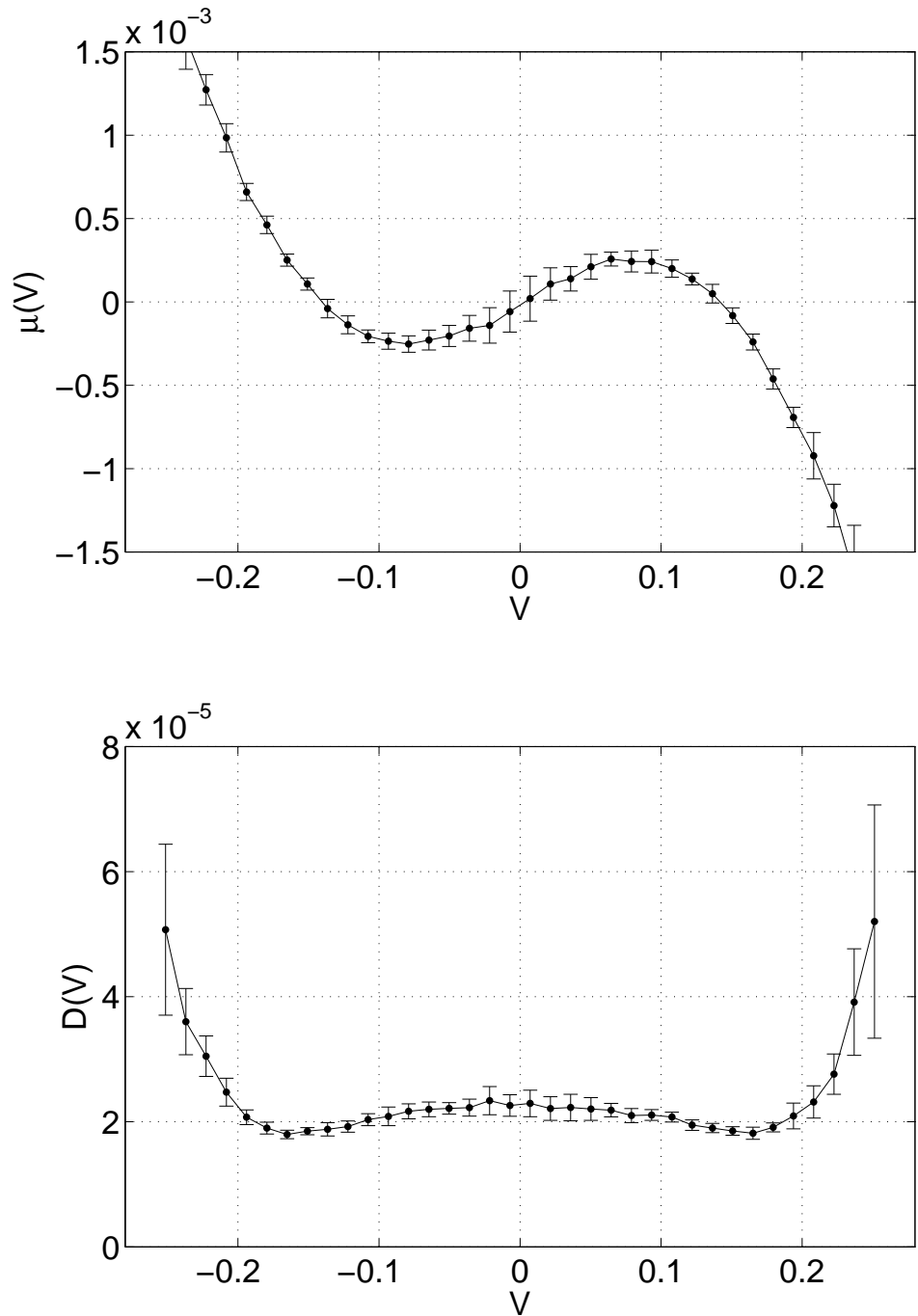}
\caption{Estimation of $\mu(V)$ and $D(V)$ by the statistics of $V$ occurrences in a long simulation.
Top: estimates of $\mu(V).$ Bottom: estimates of $D(V)$. Ten runs of length 200,000 time units
were used. Parameters are $A=0.17,\eta=10^{-4}$. Note that the error bars are smallest for
more probable values of $V$.}
\label{fig:muD}
\end{figure}

\begin{figure}
\leavevmode
\epsfxsize=5.9in
\epsfbox{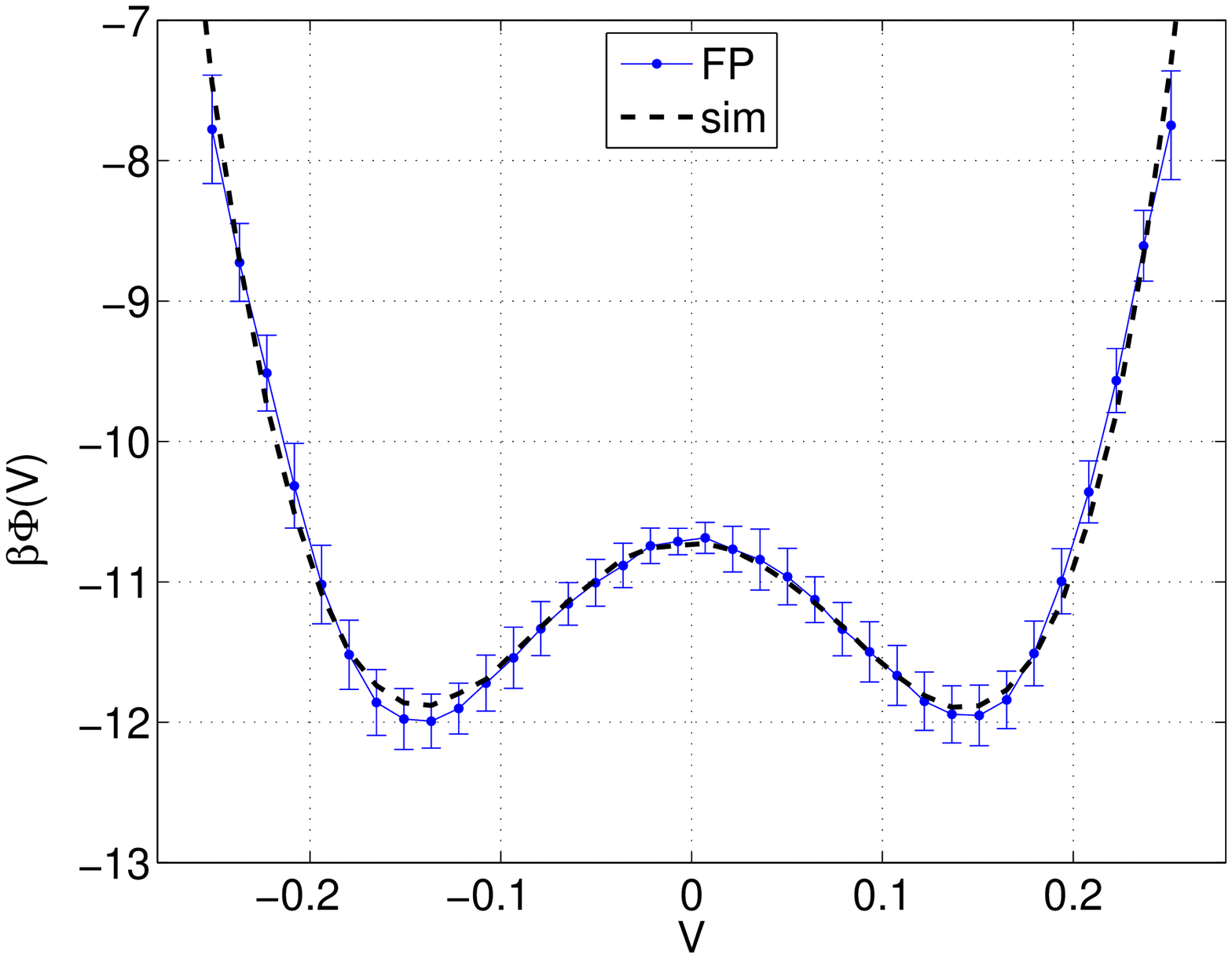}
\caption{The effective potential, $\Phi(V)$, estimated in different ways.
Dashed line from assuming that $P(V)\propto\exp{[-\beta\Phi(V)]}$; points
from assuming a Fokker-Planck equation and reconstructing $\mu(V)$ and
$D(V)$ bythe statistics of $V$ occurrences in a long simulation, then evaluating the integral~(\ref{eq:phi}).
Refer to figure~\ref{fig:muD}.}
\label{fig:pot}
\end{figure}

To establish the connection between $V$ and the effective potential $\Phi(V)$ we assume that $V(t)$ satisfies an
(unknown) Langevin equation.
Then, the probability density $P(V,t)$ satisfies a Fokker-Planck equation of the form
\begin{equation}
   \frac{\partial P(V,t)}{\partial t}=\left[-\frac{\partial}{\partial V}\mu(V)+\frac{\partial^2}{\partial V^2}D(V)\right]P(V,t).
\end{equation}
The (unknown in closed form) quantities $\mu(V)$ and $D(V)$ are estimated from ensembles of brief
simulations initialized at $V$:
\begin{equation}
   \mu(V)=\lim_{\Delta t\rightarrow 0}\langle\Delta V\rangle/\Delta t \hspace{20mm}
   2D(V)=\lim_{\Delta t\rightarrow 0}\langle[\Delta V]^2\rangle/\Delta t
\end{equation}
where $\Delta V=V(t+\Delta t)-V(t)$ and the average is taken over the realization ensemble.
More accurate, maximum likelihood based estimation techniques \cite{ait99} can
also be used.
At stationarity we have $P(V)\propto\exp{[-\beta\Phi(V)]}$ where $\Phi(V)$ is the effective
potential, related to $\mu(V)$ and $D(V)$ by~\cite{erbkev06,srikev05}
\begin{equation}
   \beta\Phi(V)=\mbox{const}-\int_0^V\frac{\mu(s)}{D(s)}ds+\log{D(V)} \label{eq:phi}.
\end{equation}
Data for the estimation of $\mu(V_0)$ and $D(V_0)$ for a given value of $V_0$ can be obtained
in two ways.
Firstly, we find occurrences of this
value in a long simulation; we then track $V(t)$ over a subsequent time interval (we use 14 time units)
for each of these occurrences, and then average the appropriate quantities.
This approach has been used previously to estimate
drift and diffusion terms in Langevin equations~\cite{grasie00}.
Typical results of estimating $\mu(V)$ and $D(V)$
using this approach
are shown in figure~\ref{fig:muD}.
Performing this for a grid of $V_0$ values provides enough data for the numerical
estimation of the integral in~(\ref{eq:phi}); the results
for the estimates in figure~\ref{fig:muD} are shown
in figure~\ref{fig:pot}.

\begin{figure}
\leavevmode
\epsfxsize=4.8in
\epsfbox{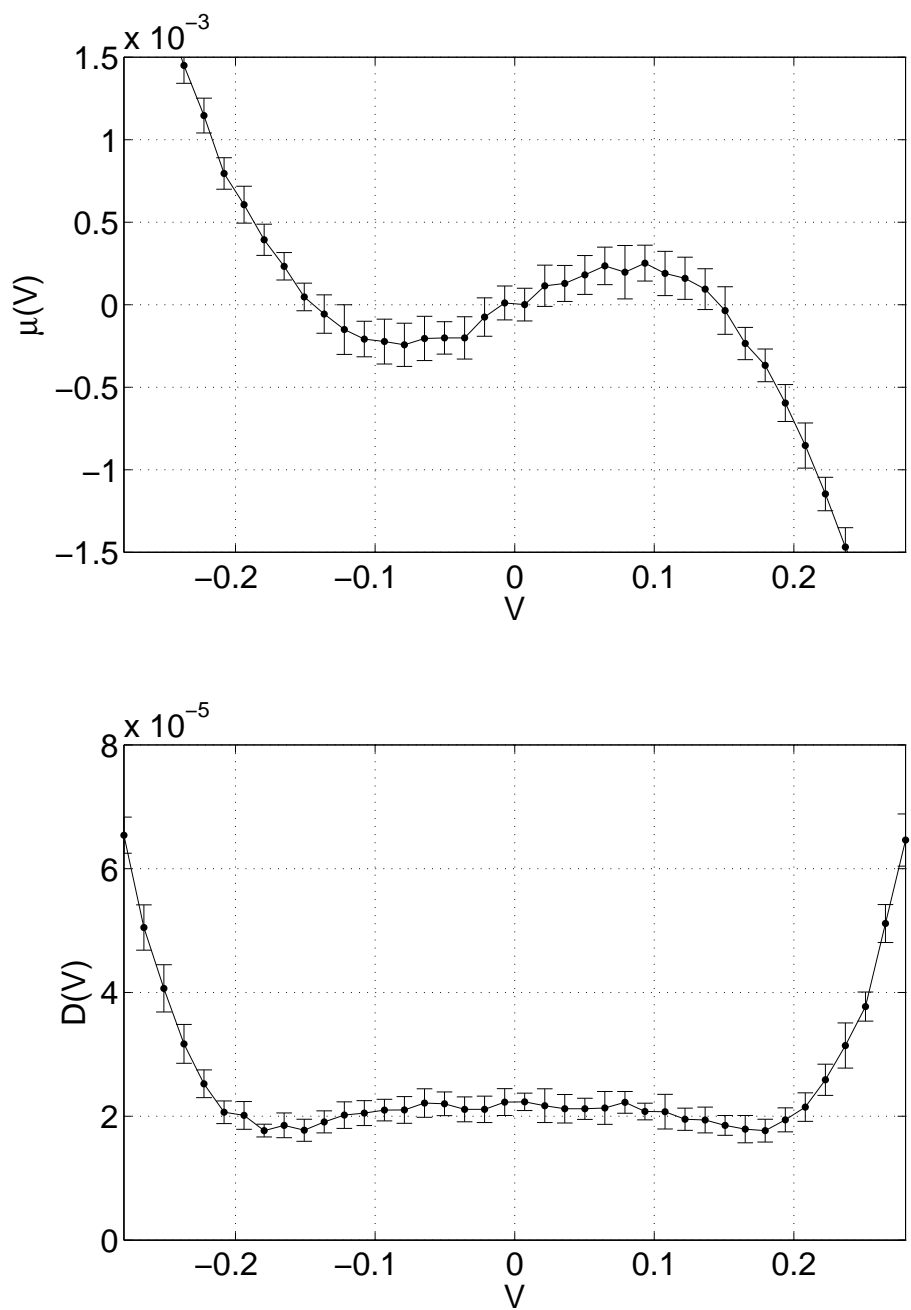}
\caption{Estimation of $\mu(V)$ and $D(V)$ by initializing the system at a particular value of $V$ and then
running for a short time. Top: estimates of $\mu(V).$ Bottom: estimates of $D(V)$.
For each value of $V,$ 4000 short bursts were run.
Parameters are $A=0.17,\eta=10^{-4}$. Note the more uniform errorbar sizes compared with those
in figure~\ref{fig:muD}.}
\label{fig:muDcold}
\end{figure}

\begin{figure}
\leavevmode
\epsfxsize=5.9in
\epsfbox{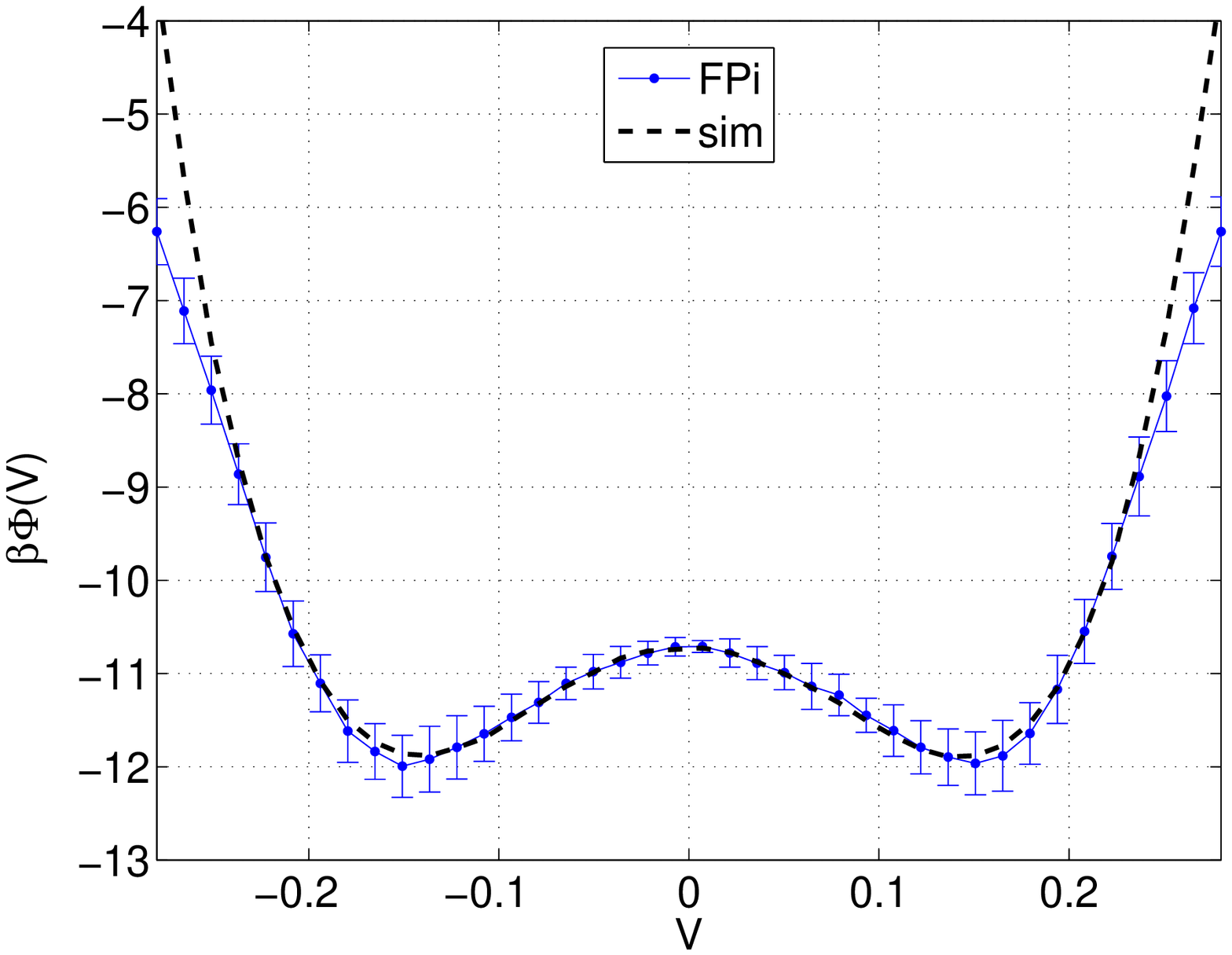}
\caption{The effective potential, $\Phi(V)$, estimated in different ways.
Dashed line from assuming that $P(V)\propto\exp{[-\beta\Phi(V)]}$, points
from assuming a Fokker-Planck equation and estimating $\mu(V)$ and $D(V)$ by initializing
the system with a particular value of $V$ and observing its short-term
behaviour. Refer to figure~\ref{fig:muDcold}.}
\label{fig:poti}
\end{figure}

In the absence of a database from a long enough (equilibrium) simulation, a
a better --- more economical --- way to find
$\mu(V)$ and $D(V)$, and thus $\Phi(V)$,
is to deliberately initialize the system with a given value of $V_0$ and then track $V(t)$
over a subsequent time interval (we again use 14 time units), repeating and averaging as above.
The advantage of this approach, as compared with the previous one, is that it need only be performed
for as many $V_0$ values as necessary for a given accuracy in the integral evaluation.
Furthermore,
for each such $V_0$ value we can reinitialize as many independent runs as necessary for accurate
estimation at will, without requiring a simulation so long that rare values
of $V_0$ are revisited enough times.
Results of the latter type of calculation are shown in figure~\ref{fig:muDcold}, and figure~\ref{fig:poti}
shows the reconstructed potential from the estimates in figure~\ref{fig:muDcold}.
The only
significant difference in this case is that initializing the system with a particular value of $V$
seems to result in a lower estimate of both $D(V)$ and $\Phi(V)$ when $|V|$ is large.
Recall, however, that these regions are rarely visited in a simulation, and thus the discrepancy
may be attributed to an insufficiently long database for the first approach.

\subsection{Bifurcations}

\begin{figure}
\leavevmode
\epsfxsize=5.9in
\epsfbox{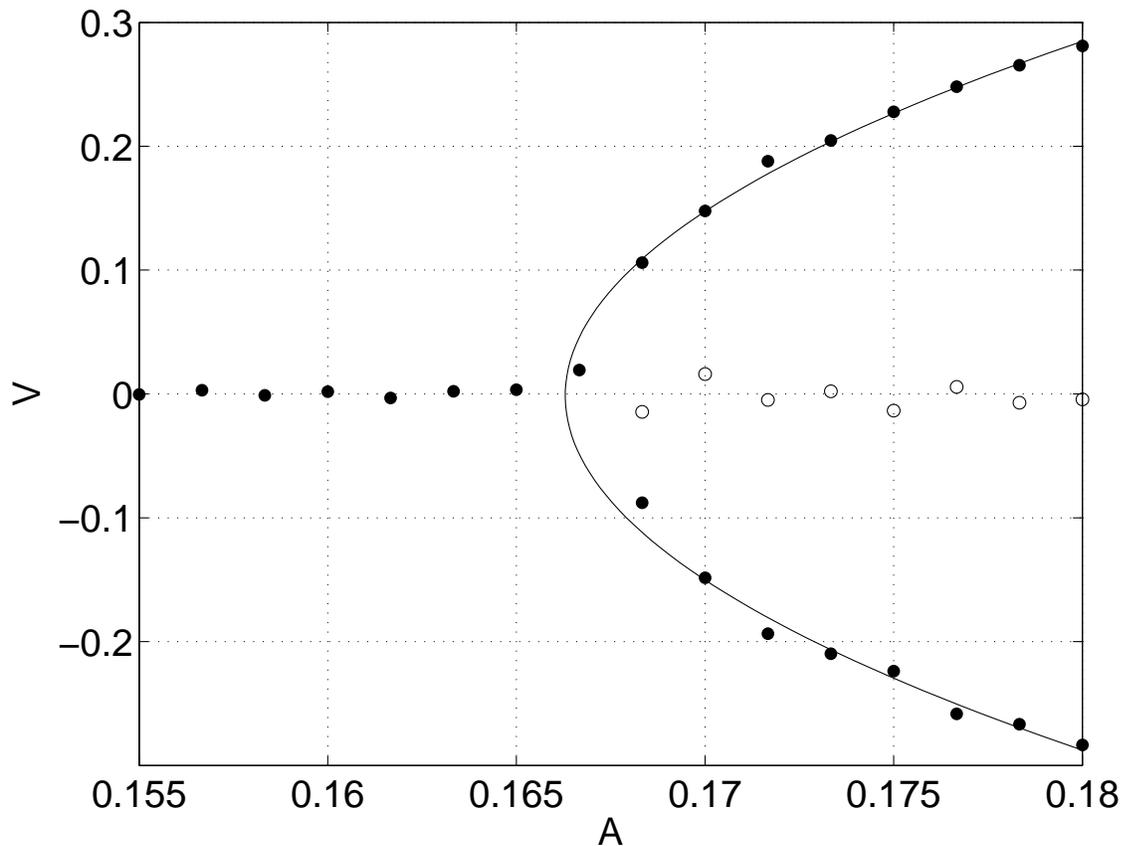}
\caption{Zeros of $\mu(V)$ as a function of $A$, for $\eta=10^{-4}$ and $M=100$. 
Filled circles indicate stable solutions; empty circles, unstable.
Also shown is a quadratic curve (in $V$) that best fits the outer branches.
For each estimation of $\mu(V)$,
10,000 short bursts were run.}
\label{fig:bifA}
\end{figure}

\begin{figure}
\leavevmode
\epsfxsize=5.9in
\epsfbox{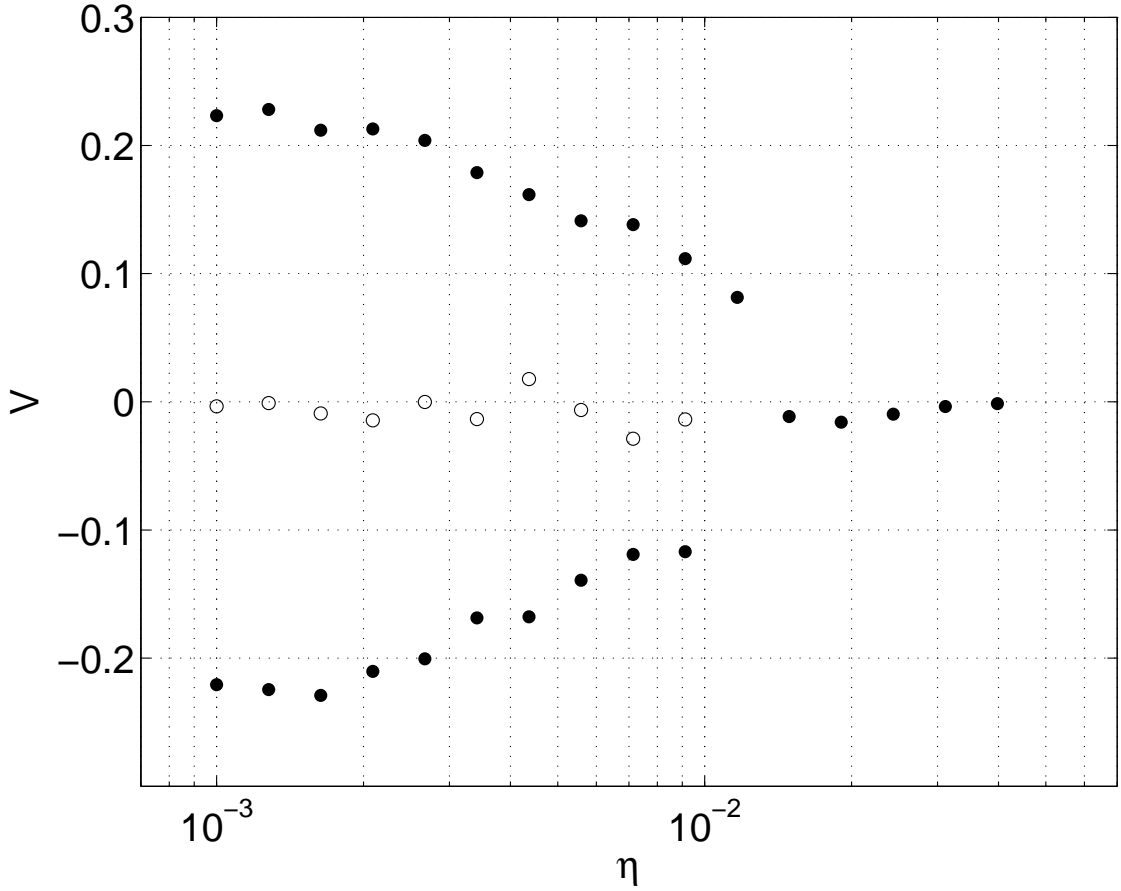}
\caption{Zeros of $\mu(V)$ as a function of $\eta$, for $A=0.175$ and $N=100$. Filled circles
indicate stable solutions; empty, unstable. For each estimation of $\mu(V)$,
10,000 short bursts were run.}
\label{fig:bifeta}
\end{figure}

A traditional bifurcation diagram for the deterministic problem with respect to a parameter such
as $A$ would involve tracing branches of steady states and constant shape travelling waves (which
can also be reduced to solutions of a fixed point problem) using pseudo-arclength and branch
switching techniques (e.g. \cite{auto}).
In the coarse-grained stochastic case, it is natural to trace instead the zeros of the drift
$\mu(V)$ as the same parameter $A$ is varied, using the same standard bifurcation codes.
It is worth noting that for the case of state-dependent noise, the local maxima of the
probability density do not exactly correspond to zeros of $\mu(V)$ and
one needs instead to find fixed points of the effective potential
given in (\ref{eq:phi}) (zeros of the right hand side of
(\ref{eq:phi}) differentiated with respect to $V$).
For the scalar case, one can implement secant-type iterative methods to converge to
the zeros of the appropriate function using function estimates only; for the multivariable
case matrix-free iterative techniques like Newton-Krylov GMRES can be used (see \cite{kel95,qiaerb06}).
In the cases we study $\mu(V)$ is well approximated by a cubic function; we take advantage of this
simplification by estimating $\mu(V)$ for just four values of $V$,
uniquely defining this cubic whose zeros we can then easily find.
Such a "noisy" bifurcation diagram is shown in figure~\ref{fig:bifA}, where we vary $A$.
Stability can be determined from the local slope, $\mu'(V)$,
at the fixed point, and is indicated in the figure.
Parameters of a different nature, such as the noise intensity or colour, can also be varied
in this context.
Figure~\ref{fig:bifeta} shows the noise-induced ``effective pitchfork" bifurcation that occurs
as $\eta$ is varied.
The location of these effective bifurcations was verified by running long simulations in their
vicinity (results not shown).
Note that if we assume in advance that the potential $\Phi(V)$ is symmetric (as the
underlying system is), the cubic function representing $\mu(V)$ would have no quadratic or constant
terms, so that only two evaluations of $\mu(V)$ would be needed to approximate the required cubic.

\subsection{Switching times}

\begin{figure}
\leavevmode
\epsfxsize=5.9in
\epsfbox{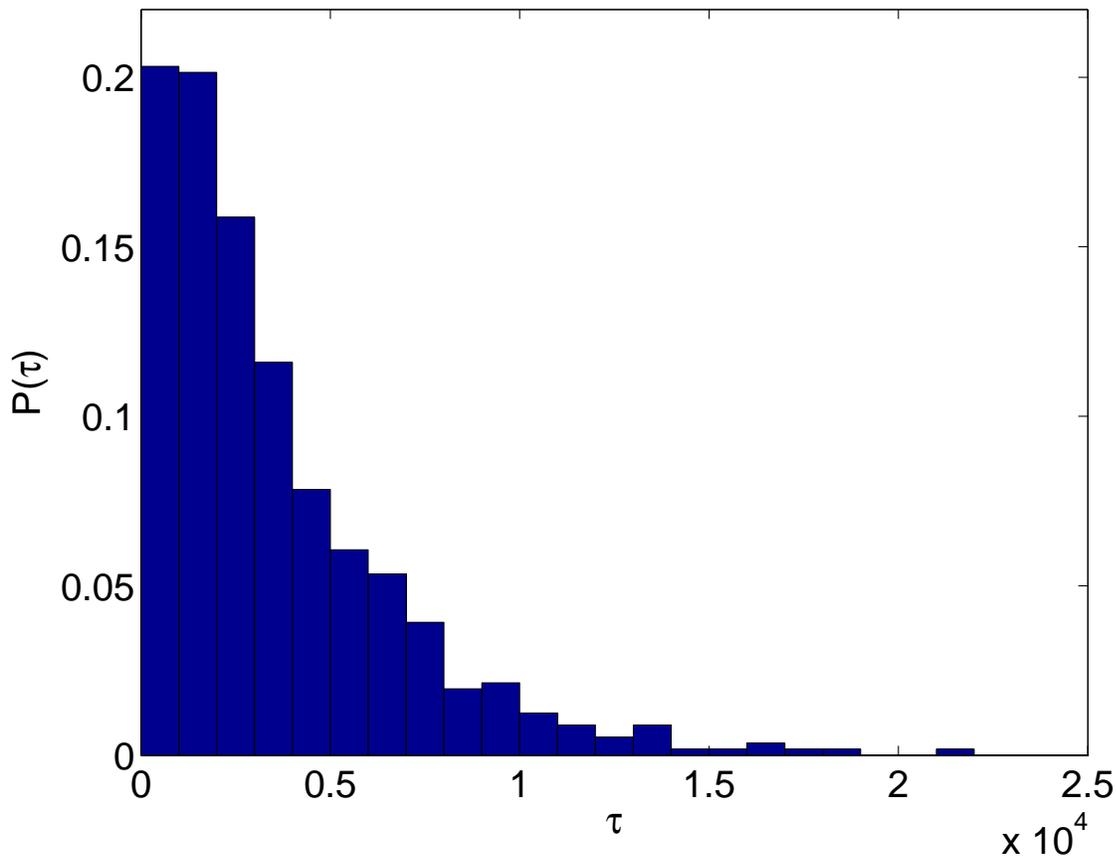}
\caption{A typical distribution of waiting times. The mean is $3.5\times 10^{3}$. Parameters
are $\eta=10^{-4}, A=0.17$.}
\label{fig:switch}
\end{figure}

One of the most important statistics of our problem is the average time between changes of
direction, or, more generally, the statistics of these switching times for the bump.
Based on our effective Fokker-Planck model and Kramers' theory,
once we have $\beta\Phi(V)$
we can also estimate the average time between switches in direction~\cite{erbkev06}, as
\begin{equation}
   \tau\approx \frac{2\pi\exp{[\beta\Delta\Phi]}}{\beta\overline{D}\sqrt{-\Phi''(V_{min})\Phi''(0)}}
\end{equation}
where $V_{min}$ is the value of $V$ at which $\Phi$ has a minimum, $\overline{D}=[D(0)+D(V_{min})]/2$
and $\Delta\Phi=\Phi(0)-\Phi(V_{min})$.
Measurements of $\overline{D}$ and
$\Phi$ for $\eta=10^{-4}$ and $A=0.17$ give $\tau\approx 3\times 10^3$.
A typical distribution of waiting times from a long
simulation is shown
in figure~\ref{fig:switch}.
The mean for the data shown is $3.5\times 10^{3}$, in excellent agreement with our Kramers' approximation.

\begin{figure}
\begin{tabular}{c}
(a) \\
\epsfxsize=3.9in
\centerline{
\epsfbox{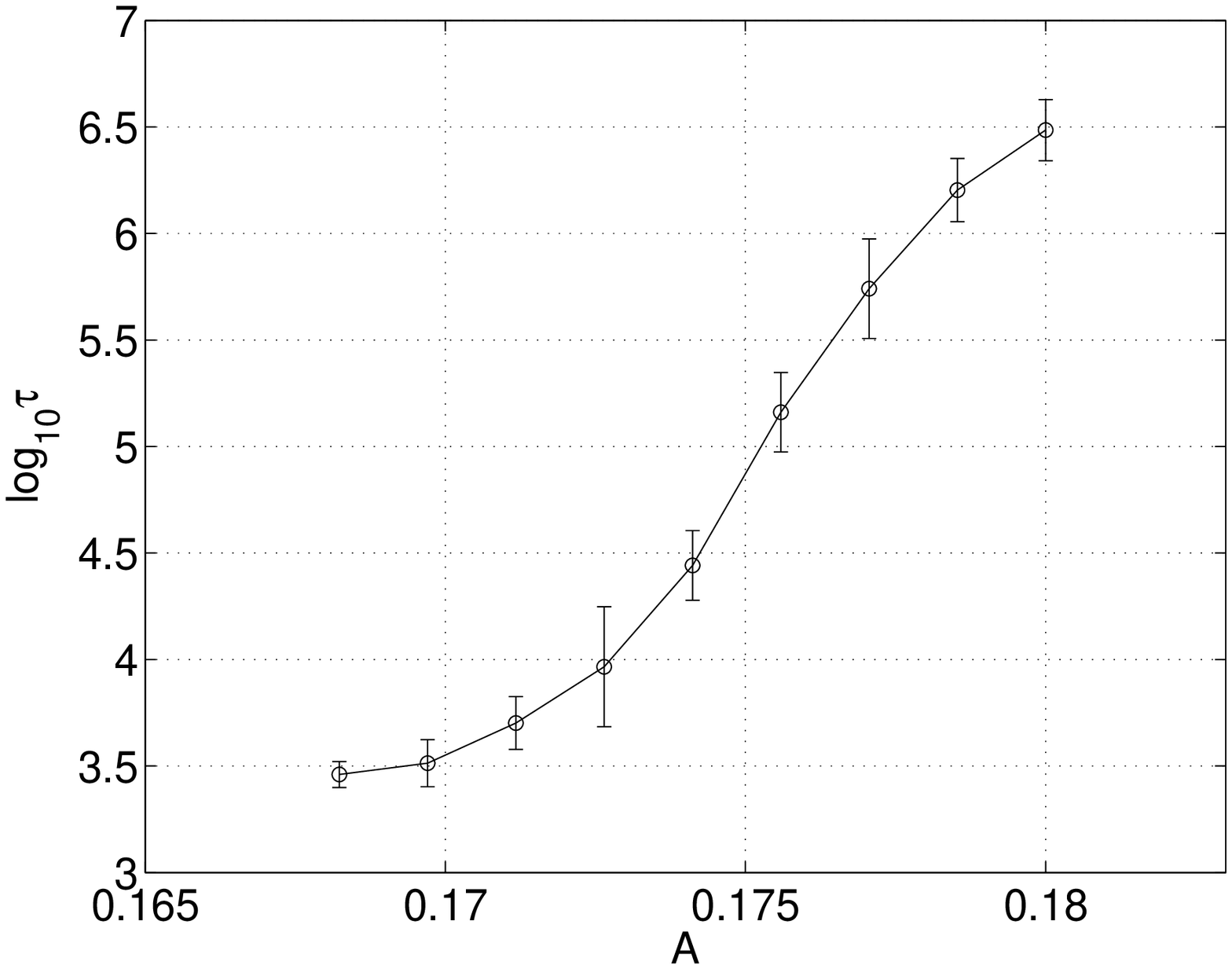}} \\
(b) \\
\epsfxsize=3.9in
\epsfbox{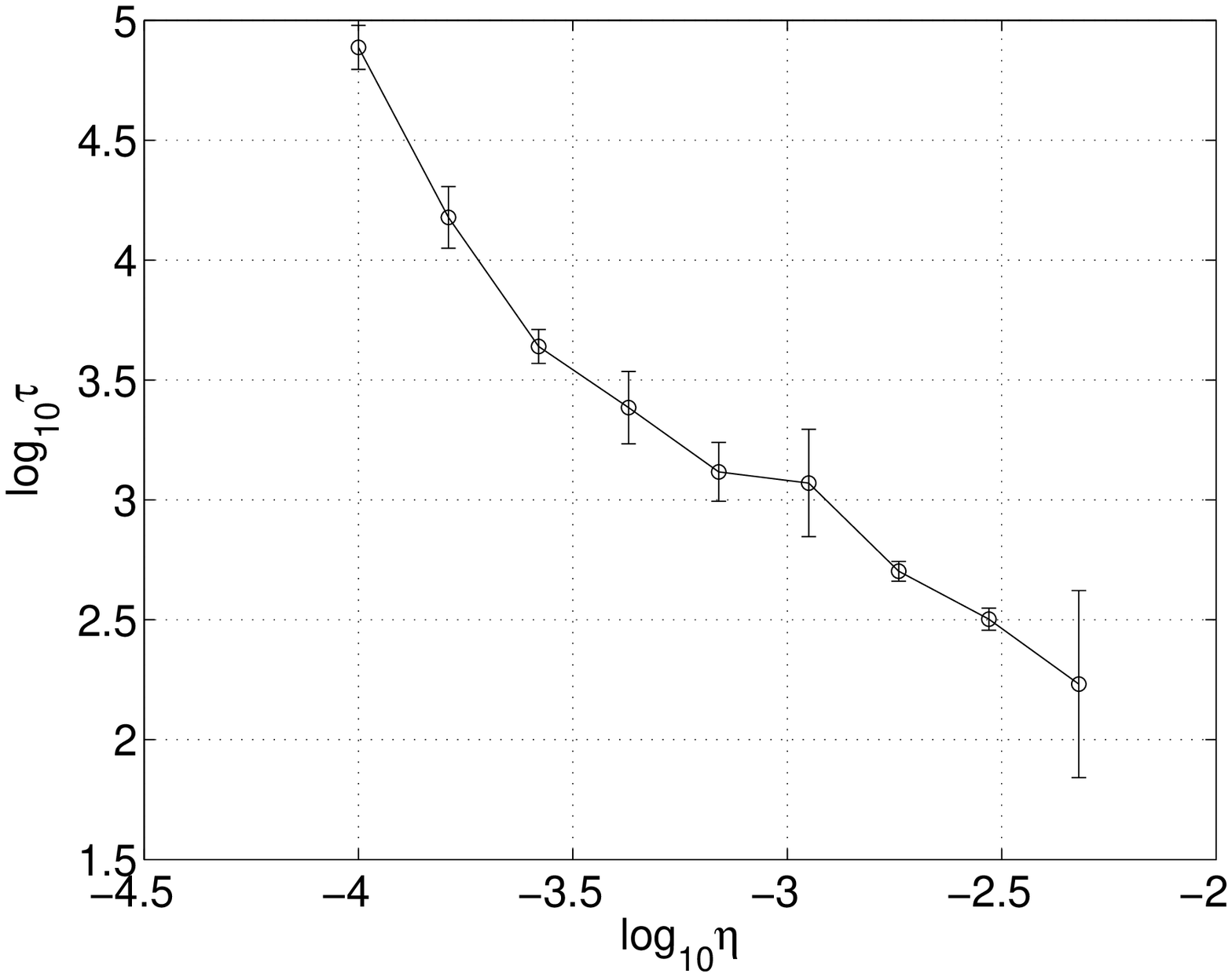}
\end{tabular}
\caption{(a) Switching time, $\tau$, plotted on a log scale, as a function of $A$. $\eta=10^{-4}.$ (b)
Switching time, $\tau$, plotted on a log scale, as a function of $\eta$ (also on a log scale).
$A=0.175.$
For both panels, $\tau$ was estimated
from a calculation of $\beta\Phi$, which itself was computed by initializing $V$ at appropriate values
and then running to observe its short-term evolution.}
\label{fig:tauA}
\end{figure}

Along the same lines, we can quantify how the average switching time, $\tau$, depends on parameters.
Typical results are shown in figure~\ref{fig:tauA}.
%
Here, the effective potential is not estimated from a long simulation database, but rather ``on demand"
by initializing the system at particular values of $V$, performing short time simulation bursts, and
processing their results as before.
A clear advantage of this approach, as opposed to running the system
for long enough to measure sufficient occurrences of
switches, is that large values of $\tau$ can be inferred from a reasonable number of short simulations.
See
figure~\ref{fig:tauA}, for example, where we can ``measure'' switching times greater than $10^6$.

For a simulation of a given, fixed duration $T$,
these curves can be used to infer the parameter values for which
we expect to see direction changes.
If $\tau>T$, we will
typically see the bump move in one direction only, whereas if $\tau<T$ we will see switching,
leading to a drop in the absolute value of the distance travelled during the time $T$.
For example,
from figure~\ref{fig:tauA} (b) we see that for simulations of length $T=10^3$, $\eta$ must be greater
than approximately $10^{-3}$ to observe any switching of direction.


\subsection{White versus coloured noise}
In~\cite{lailon01} the authors investigated the effects of adding coloured noise
to~(\ref{eq:dudtA})-(\ref{eq:dadtA}).
They found that increasing the correlation time
of the coloured noise made the bump less likely to move.
This was rationalized in terms of ``frozen'' noise
(i.e.~spatial inhomogeneity) which was likely to ``pin'' the bump, preventing it from moving.
Increasing
the correlation time of the coloured noise was thought of as interpolating between
Gaussian white noise (with delta function autocorrelation in time) to frozen, spatially structured, noise.
Here we investigate the effect
again, quantifying the influence of noise colour on the location of the underlying bifurcation.

\begin{figure}
\leavevmode
\epsfxsize=5.9in
\epsfbox{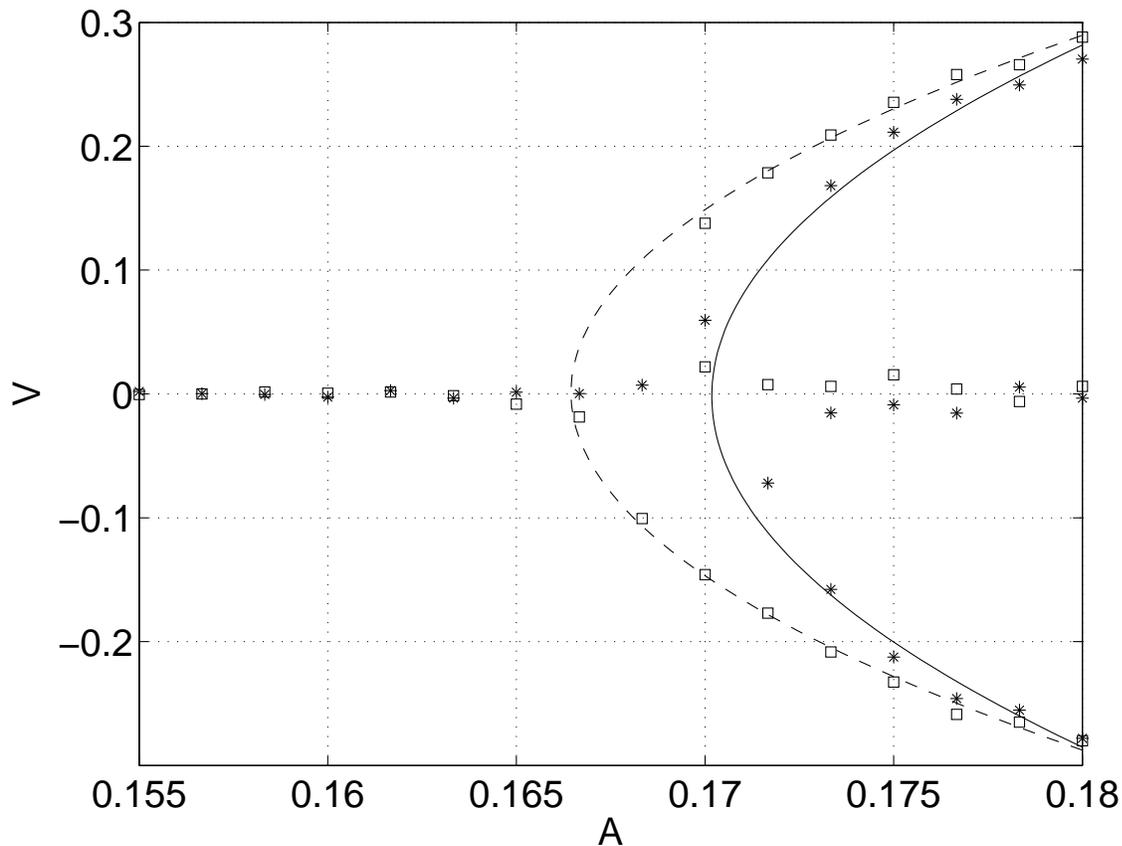}
\caption{Zeros of $\mu(V)$ as a function of $A$ for coloured noise with correlation time
$\lambda=100$ (solid line and asterisks) and $\lambda=1$ (dashed line and squares).
The curves are quadratic functions in $V$ fitted to the data points.}
\label{fig:color}
\end{figure}

As in~\cite{lailon01} we add noise to the system by adding a term
$\eta_i(t)$ to each of~(\ref{eq:dudtA}), with $\langle\eta_i(t)\rangle=0$ and
$\{\langle\eta_i(t)\eta_j(s)\rangle\}=2\epsilon\nu_{ij}e^{-|t-s|/\lambda}$, where $\nu_{ij}=0$
if $i\neq j$ and
1 if $i=j$.
The notation $\{\cdots\}$ indicates averaging over the initial distribution of $\eta(0)$
values, taken from a Gaussian with zero mean and variance $2\epsilon$.
We keep $\epsilon=10^{-4}$
and vary the correlation time $\lambda$.
The results are shown in figure~\ref{fig:color} where we have
plotted zeros of $\mu(V)$ as a function of $A$ for $\lambda=1$ and $\lambda=100$.
The figure was computed
using the short runs initialized at prescribed $V$ values.
The results are consistent
with those found in~\cite{lailon01}, where it was seen that coloured noise was more effective at
slowing the bump than Gaussian white noise of the same power, and that the effect was
stronger for longer correlation times.
The results presented here show that
this effect can be rationalized in terms of shifting a bifurcation point.
Results for several
other values of $\lambda$
indicate that as $\lambda$ is increased, the bifurcation moves to higher values of $A$ (results not shown).

\section{Diffusion maps and the data-based detection of coarse observables}

\begin{figure}
\leavevmode
\epsfxsize=3.5in
\centerline{
\epsfbox{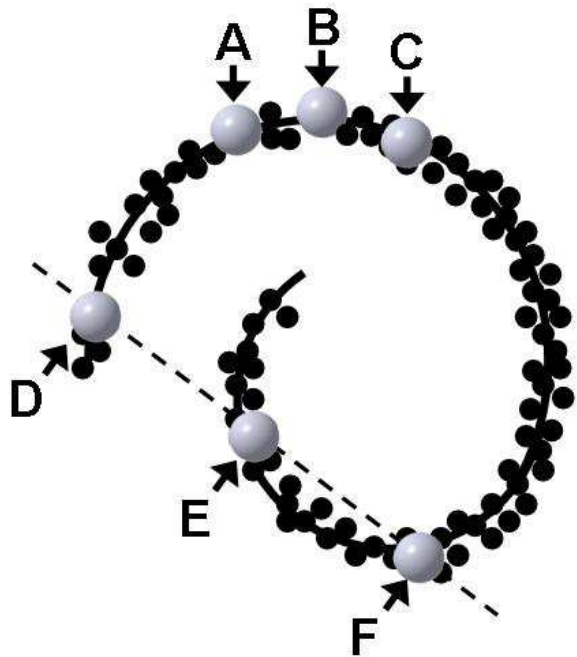}
}
\caption{Schematic of datapoints (filled circles) in
  $\mathbb{R}^{2}$ that lie along a curve (solid line). Euclidean
  distance is a good measure of the separation between points A, B,
  and C along the curve but is a poor measure of the separation of
  points D, E, and F along the curve.}
\label{fig:dataswirl}
\end{figure}

The previous results relied on our ability to choose a scalar variable, $V$, whose
value correlates with the state of the high-dimensional system (after brief initial transients
have equilibrated).
However, we are often
faced with dynamical systems for which choosing such a variable is far from obvious.
We then need a systematic procedure for determining, from the results of a simulation,
a good low-dimensional representation of the system state.
Such data-mining, so-called ``manifold learning" techniques have been a focus of
intense research in recent years \cite{coilaf05A,coilaf05B,erbfre06}; we will use here the recently
developed diffusion map approach~\cite{nadlaf05,nadlaf06}.
A large data ensemble from direct simulations of our system can be represented as
a cloud of points in a high-dimensional space (here,
$\mathbb{R}^{200}$); see the schematic in figure~\ref{fig:dataswirl}.

When two such points (activity profiles) are very close to each other (i.e. when their
Euclidean distance is small, as in the case of points A, B and C in the schematic),
we can consider this distance as representative of the
``intrinsic similarity" of the two configurations - in some sense, of how easy it is
for the system dynamics to cause a transition from one configuration to the other.
When, however, this Euclidean distance is larger than some threshold (as in the case of
points D, E and F in the schematic) the Euclidean distance stops being
a good measure of the ``intrinsic similarity" between configurations -- the ``effort to
transition" from E to F, measured by the arclength between them, is
clearly less than the effort to transition from E to D, even though the Euclidean distances
DE and EF are similar.
The main idea underpinning diffusion maps is to perform a random walk on a graph
in which data points are vertices, and connection strengths between data points are
given by a Gaussian kernel of the form $K\left({\bf x},{\bf y}\right)=\exp\left(-\norm{{\bf x}-{\bf y}}^2/\sigma^2\right)$.

When the points are farther away than a cutoff (controlled by the
parameter $\sigma$)
the vertices are effectively disconnected; when they are very close, the strength
of the connection is controlled by their (small) Euclidean distance.
A random walk on such a graph gives rise to a {\em diffusion distance} --
starting from the same source point and diffusing on the graph for some time,
we then look for equal density contours; points on such a contour
are equally easy to access from the source point, and are therefore at equal
{\em diffusion distance} from it, even though their Euclidean distances from
the source may vary substantially.
The procedure also provides a set of
transformed coordinates; the Euclidean distance {\em in these new coordinates}
is a true measure of intrinsic datapoint similarity.
The procedure can be thought of as a nonlinear generalization of Principal
Component Analysis \cite{jol86}.

\begin{figure}[t]
\centerline{
\psfig{file=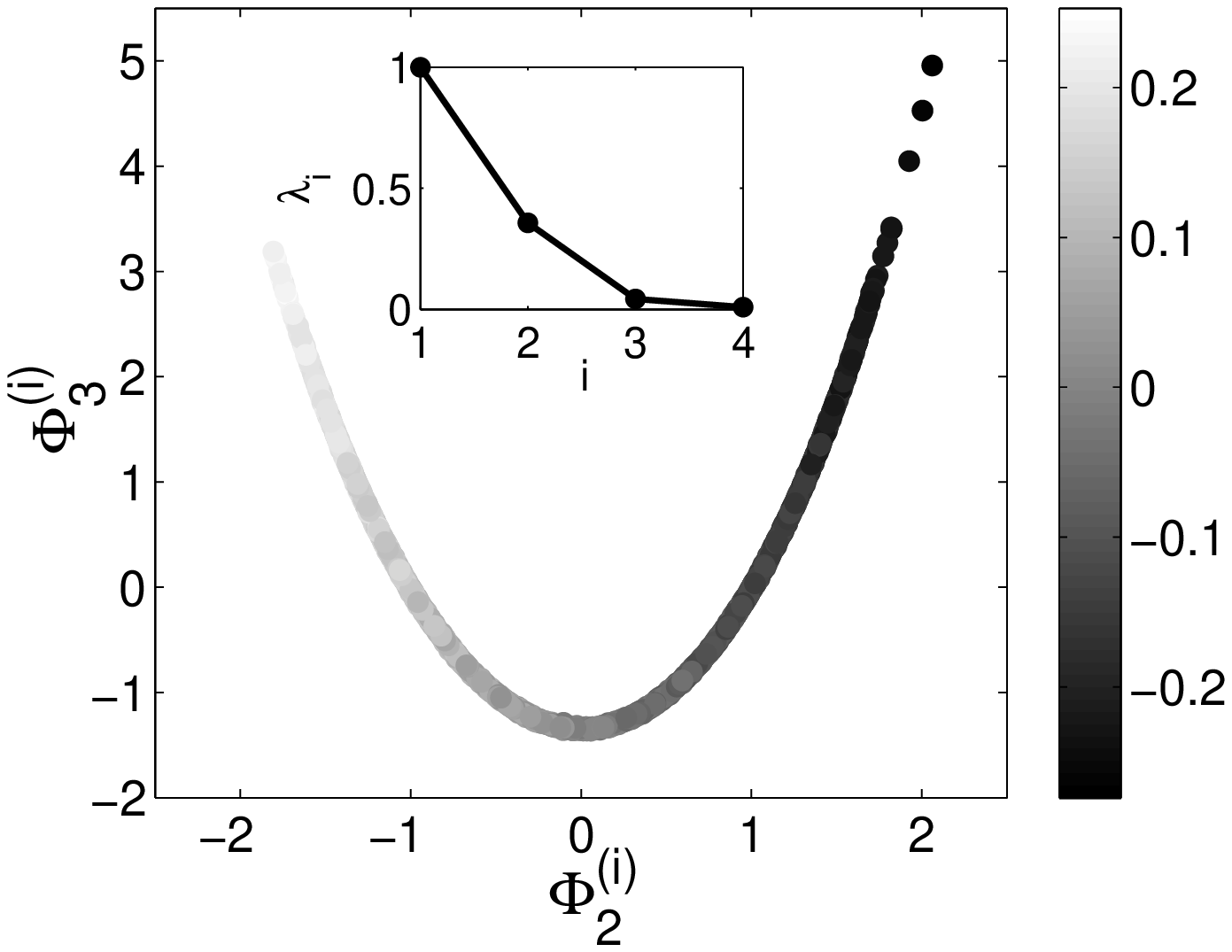,height=2.5in}
\psfig{file=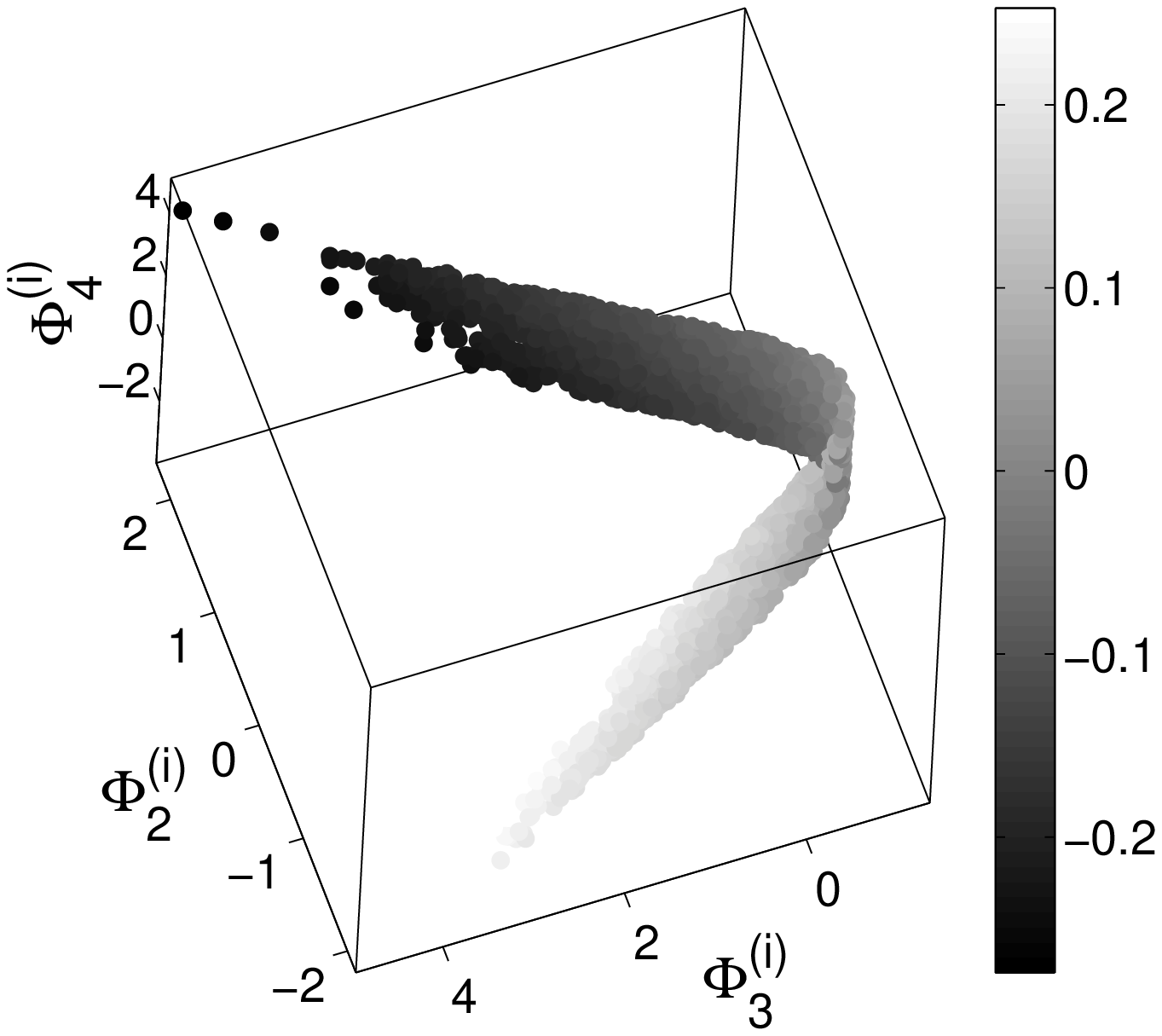,height=2.5in}
}
\caption{ Left panel: Diffusion map plotting components of the top
    two significant eigenvectors of Markov matrix ${\bf M}$ (defined
    in Appendix)
    constructed from simulation data. Datapoints
    are shaded by (and appear ordered according to) the value of the
    empirical coordinate $V(t)$ (defined in
    Section \ref{sect:intro}). The inset shows the leading eigenvalue spectrum
    for the diffusion map. Right panel: 3D diffusion map plotting components of top
    three significant eigenvectors
}
\label{fig:2+3ddmap}
\end{figure}

\begin{figure}
\leavevmode
\epsfxsize=5.in
\centerline{
\epsfbox{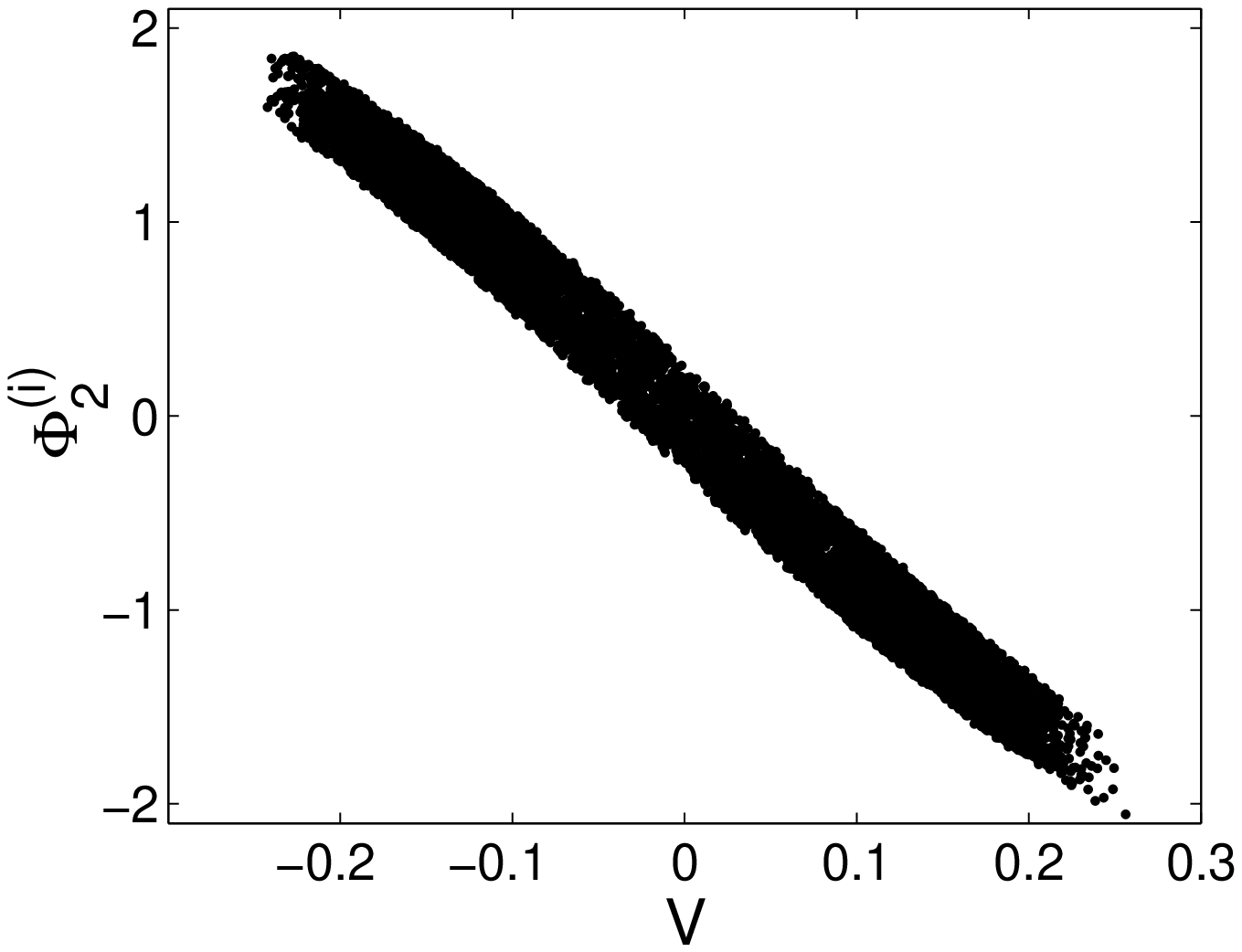}
}
\caption{Coordinate $V(t)$ versus top significant eigenvector
  coordinate $\Phi_2^{(i)}$ for long simulation results.}
\label{fig:coordMAP}
\end{figure}

In our case, we start from a long simulation (of 30,000 time units) during
which we sample the $u(x),a(x)$ profile every 8 time units,
giving $N=3,\!750$ data points. (From now on, we use the same parameter values as those in
figure~\ref{fig:example}.)
Using the procedure briefly outlined in the Appendix we process these data
forming a Markov matrix whose leading eigenvectors provide useful ``reduction
coordinates" for our data -- in particular, the second eigenvector $\Phi_2$ of
this matrix is our computer-assisted candidate coarse observable (the first
eigenvector, corresponding to the eigenvalue 1, is trivial).
Figure~\ref{fig:2+3ddmap} shows the simulation data projected on the
first two ``reduction coordinates" -- the first two nontrivial eigenvectors $\Phi_2$ and $\Phi_3$ of
the appropriate Markov matrix. The right panel of
figure~\ref{fig:2+3ddmap} shows a three-dimensional diffusion map, in
terms of the first three nontrivial reduction coordinates. The coordinates of the $i^{th}$ datapoint in this
map are $(\Phi_2^{(i)},\Phi_3^{(i)},\Phi_4^{(i)})$ . Simulation points
in the diffusion map are coloured according to their associated value
of the empirical coordinate $V$ -- it is clear that points are ordered
along the curve by their values of this known reaction coordinate. This sorting of
high-dimensional data vectors is obtained in an ``automated'' fashion
by the diffusion map calculation.
\begin{figure}
\centerline{
\epsfxsize=5.5in
\epsfbox{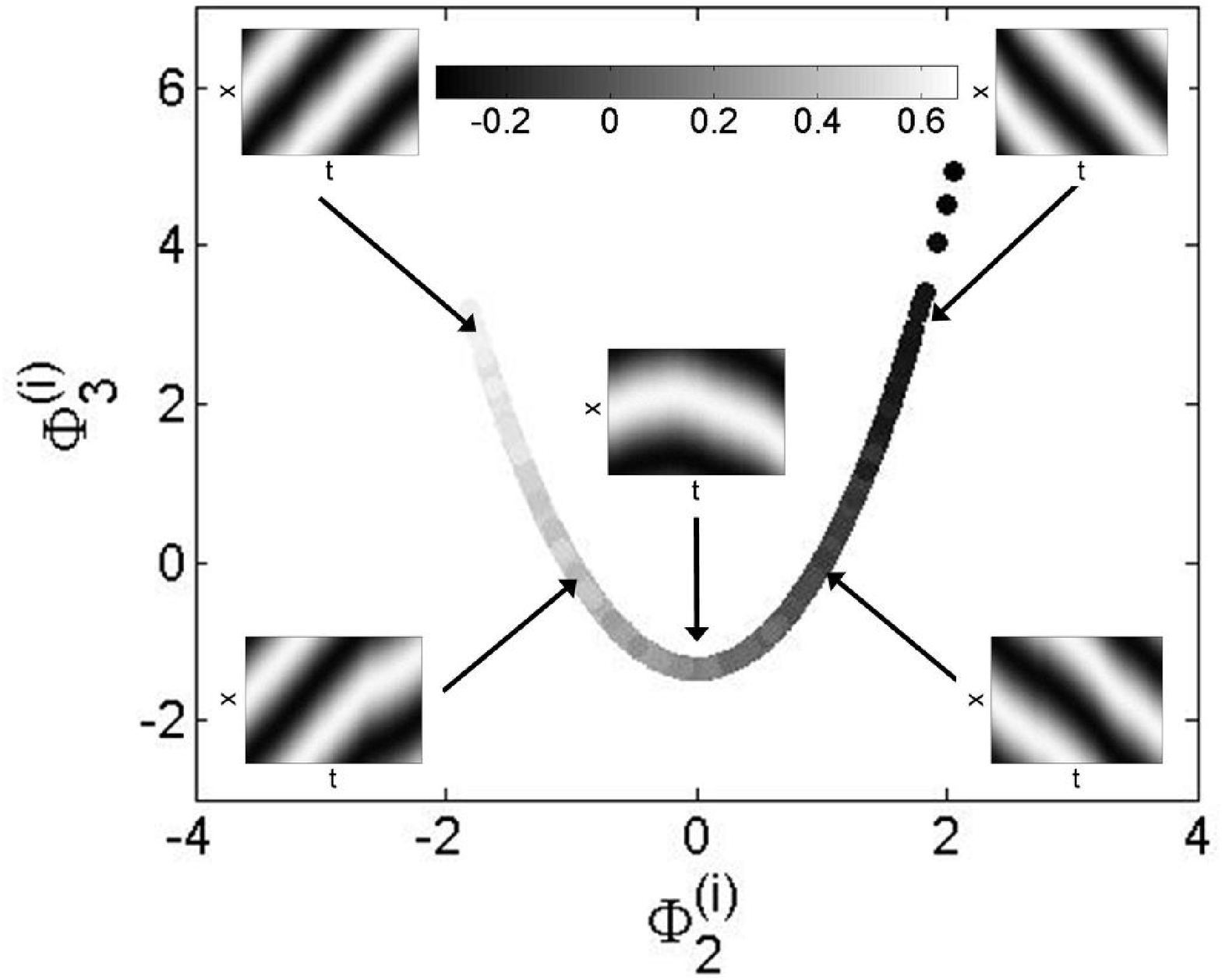}
}
\caption{ 2-dimensional diffusion map representation of long
    simulation data (as shown in
    figure~\ref{fig:2+3ddmap}). The insets show representative local  solution evolution for $u$
    profiles initialized at 5 different regions of the map
    ($\Phi_2=-1.8,-1,0,1,$ and $1.8$). Shadebar shown at top of figure
    indicates $u$ values in space-time plots.
}
\label{fig:2ddmap+evol}
\end{figure}
All data instances clearly collapse on a one-dimensional
curve. Figure~\ref{fig:2ddmap+evol} indicates five representative data
points (at values of $\Phi_2=-1.8,-1,0,1,$ and $1.8$) on this curve. 
The insets, showing brief space-time plots initialized at the corresponding
data points, clearly indicate that right-moving states reside on one end of
this curve, left moving states on the opposite end, and transition states between
the two reside in the interior of the curve.
It is clear that the curve is one-to-one with the variable $\Phi_2$
(for every value of $\Phi_2$ there exists one point on the curve); this
suggests that $\Phi_2$ would be a good scalar observable for our system.
Figure~\ref{fig:resmulti} (third row) shows the evolution of this data-based observable during our
system dynamics; it is clearly capable of capturing the direction switches
in the system in a manner comparable with the ``experience-based" coordinate
$V$. Another possible scalar observable we considered is the arc-length distance along
the curve connecting components of the top two eigenvectors. The variation in this
observable is comparable to that of $\Phi_2$ and also accurately describes the
direction of travel of the bump during the simulation (results not shown).
Figure~\ref{fig:coordMAP} plots the simulation data in terms of the
two coarse coordinates: the experience-based $V$ and the data-based $\Phi_2$.
At first sight, a one-to-one correspondence between the two observables
is apparent, suggesting that our experience-based coordinate was indeed a good
choice.

\begin{figure}
\leavevmode
\epsfxsize=5.in
\centerline{
\epsfbox{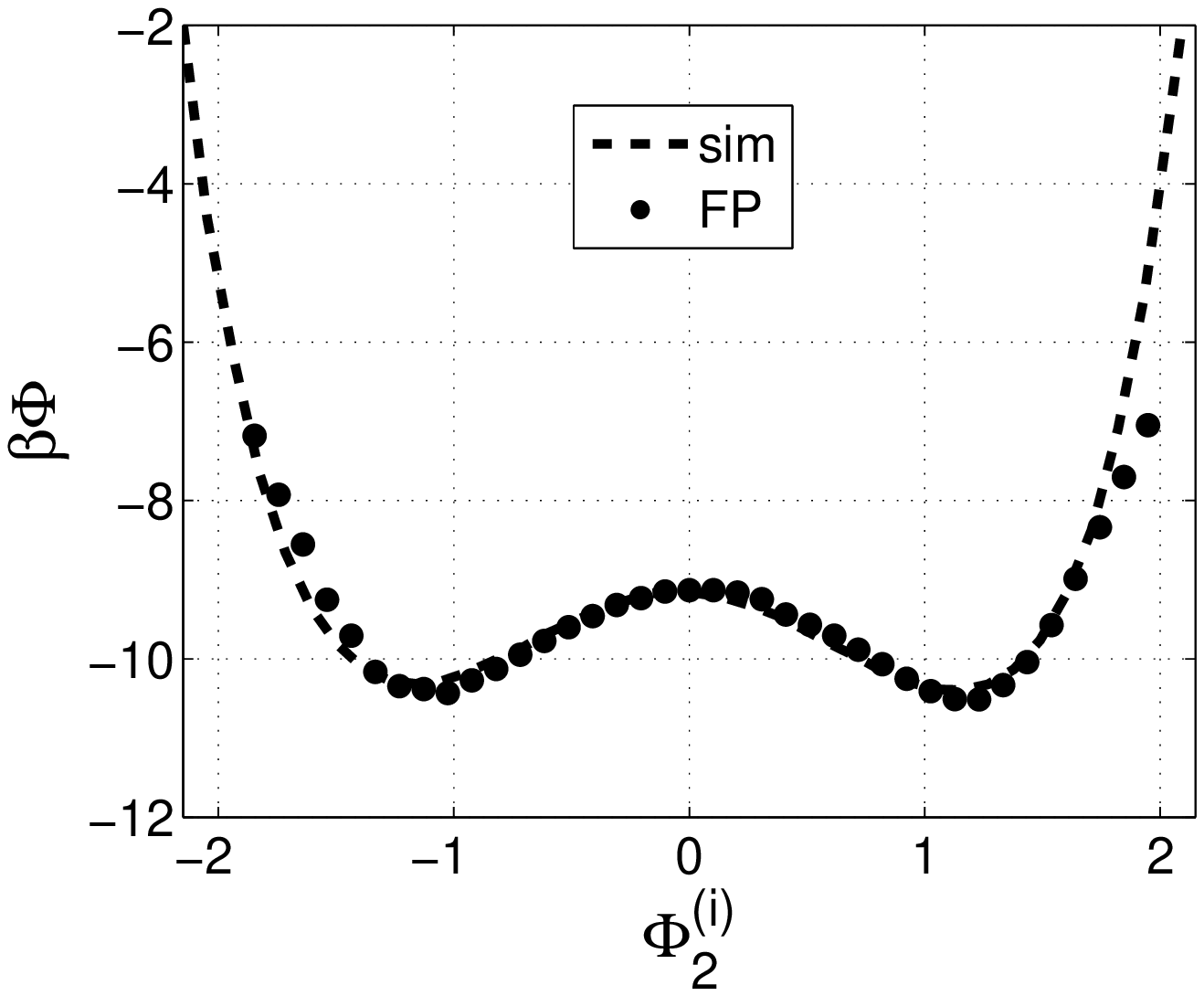}
}
\caption{Effective potential $\beta\Phi$ estimated in terms of diffusion
  map coordinate $\Phi_2^{(i)}$. Dashed line: from assuming that
  $P(\Phi_2^{(i)})\propto\exp{[-\beta\Phi]}$; points from assuming a
  Fokker-Planck equation, reconstructing $\mu(\Phi_2^{(i)})$ and
  $D(\Phi_2^{(i)})$ and computing the integral~(\ref{eq:phi}).}
\label{fig:dmapFE}
\end{figure}

The processing of the database from a long, equilibrium run simulation that gave
us an effective potential in terms of the observable $V$
can now also be performed in terms of the computer-assisted observable $\Phi_2$;
the resulting two-well effective potential and the corresponding effective Langevin
description are practically isomorphic to the
$V$-based description.
The effective potential computed in terms of the computer-assisted
observable $\Phi_2$ is shown in figure~\ref{fig:dmapFE}.

\subsection{Computing with diffusion map coordinates}

When we attempt to {\it directly} construct the effective potential in terms of the
diffusion map coordinate $\Phi_2$ (i.e., not from processing a long run database),
it becomes necessary to initialize short simulation bursts consistent with
particular values of $\Phi_2$ as well as observe the $\Phi_2$ values
corresponding to results of detailed simulation.
The diffusion map calculation that leads to the identification of the coordinate $\Phi_2$
automatically provides its value on each of the data points used in the calculation
(see the Appendix).
It therefore becomes important to establish computational protocols that routinely allow
the translation between physical coordinates and diffusion map coordinates.
This translation must be possible in both directions: to initialize a short system run
we need to translate a value of $\Phi_2$ to system initial conditions (living in  $\mathbb{R}^{200}$);
processing of the short simulation bursts involves extracting the $\Phi_2$ values of the resulting states.
These translation operations are termed ``lifting" and ``restriction" in the equation-free
framework \cite{geakev02,kevgea03}.
We now briefly discuss possible implementations of such protocols, starting with the
{\em restriction}: obtaining the $\Phi_2$ values of new system configurations, not in
the original database, that result from our short simulation bursts.

For points {\em in the original dataset} (used to assemble the neighbourhood
matrix ${\bf K}$) eigenvalues and eigenvectors of the symmetric kernel ${\bf M_s}$
are related by
\begin{equation}
 {\bf M_s \Psi}_j=\lambda_j {\bf\Psi}_j.
\label{eigenrel1}
\end{equation}
For a particular point ${\bf x}_i$ in the dataset we have 
\begin{equation}
 \sum_{k=1}^N {\bf M_s}({\bf x}_i,{\bf x}_k) \Psi_j^{(k)}=\lambda_j \Psi_j^{(i)}.
\label{eigenrel2}
\end{equation}
For a new datapoint ${\bf x}^{(new)}$ the same formula should
be valid. However, while the right hand side is unknown, the left hand
side may be computed for the new datapoint, which gives the following
formula for the extension of the diffusion map 
\begin{equation}
 {\Psi}_j^{(new)}=\frac{1}{\lambda_j}\sum_{k=1}^N {\bf
 \widetilde{M}_s}({\bf x}^{(new)},{\bf x}_k) \Psi_j^{(k)},
\label{Nystrom}
\end{equation}
where ${\bf \widetilde{M}_s}$ is a generalized kernel which extends to
new sample points and is given by
\begin{equation}
  {\bf \widetilde{M}_s}({\bf x}^{(new)},{\bf x}_k)=\frac{1}{N}
\frac{K({\bf x}^{(new)},{\bf x}_k)}
{\sqrt{E_{\bf y}\left(K({\bf x}^{(new)},{\bf y})\right)
E_{{\bf y}^{\prime}}\left(K({\bf x}_k,{\bf y}^{\prime})\right)}}.
\label{kernelGen}
\end{equation}
Equation (\ref{Nystrom}) is called the Nystr\"{o}m formula~\cite{bak77} and allows eigenvector
components associated with new data vectors ${\bf x}^{(new)}$, outside of
the original sample, to be computed by eigenspace interpolation
instead of repeated eigendecomposition of the Markov matrix augmented by
the new data vectors. 
The expectations appearing in the denominator of~(\ref{kernelGen}) are computed from the
empirical data using
\begin{equation}
 E_{{\bf y}}\left(K({\bf x},{\bf y})\right)=\frac{1}{N}\sum_{i=1}^N
K({\bf x},{\bf y}_i).
\label{expectation}
\end{equation}
Use of the Nystr\"{o}m formula (\ref{Nystrom}) is central to the efficient computation of
 diffusion map coordinates associated with a high-dimensional vector
 ({\em restriction}) and also the preparation of
 a data vector {\bf x} ({\em lifting}) with desired diffusion map
 coordinate values.
 Equation~(\ref{Nystrom}) is used to first compute the eigenvector components
 ${\Psi}_j^{(new)}$ (associated with the symmetric matrix ${\bf
 M_s}$); diffusion map coordinates (associated with ${\bf M}$) for the new data point
 may then be computed using
\begin{equation}
 \Phi_j^{(new)}=\frac{\Psi_j^{(new)}}{\Psi_1^{(new)}}.
\label{eigencompNEW}
\end{equation}
%

The Nystr\"{o}m formula allows calculation of the diffusion map eigenvector
components $\Phi_j^{new}$ associated with a new data vector ${\bf x}^{(new)}$.
A full eigendecomposition is typically performed first for a
representative subset of simulation datapoints, identifying the relevant eigenvectors.
The Nystr\"{o}m
formula is then used to perform the restriction operation in~(\ref{Nystrom}) which amounts to
interpolation in diffusion map space.

\begin{figure}
\epsfxsize=4.75in
\centerline{
\epsfbox{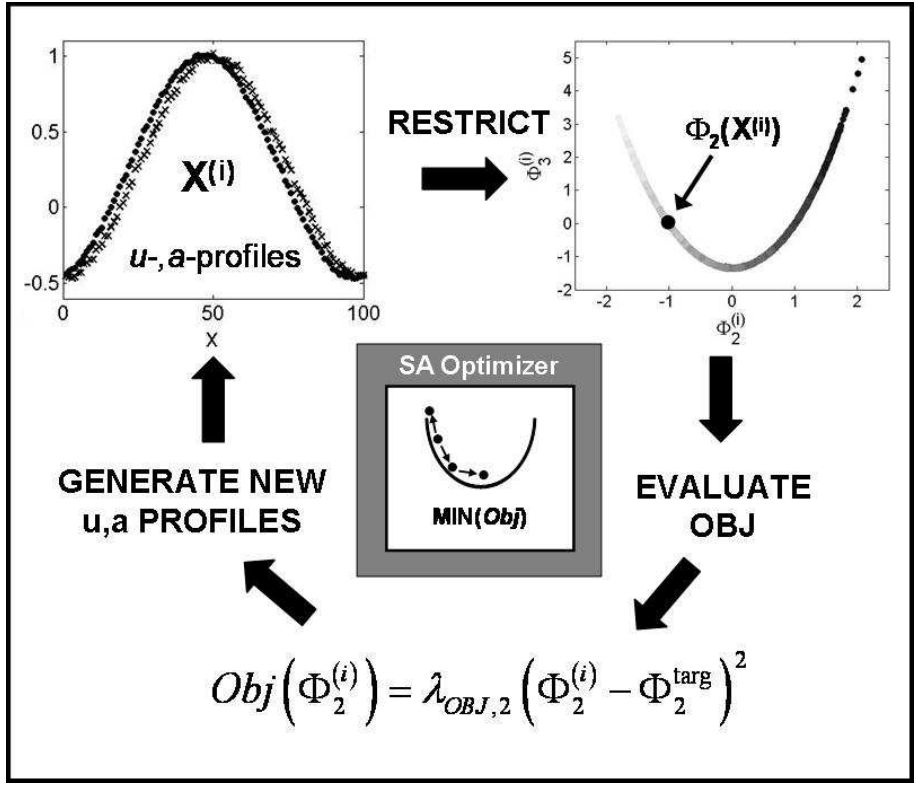}
}
\caption{ Protocol for lifting from diffusion map coordinates
    using Simulated Annealing (SA). We show (top left) normalized $u$ and
$a$ solution profiles which define the trial vector ${\bf
  x}^{(trial)}$. This 200-dimensional vector is restricted, 
the ``location'' of the data vector in the diffusion map is
identified (top right) and the corresponding objective value computed
(bottom). The objective value guides the selection of a new search
    direction (in $\mathbb{R}^{200}$), generating a new set of
    solution profiles (a new ${\bf
  x}^{(trial)}$) and the above procedure is repeated.
}
\label{figLIFTseq}
\end{figure}

Our protocol for {\it lifting} from prescribed diffusion map coordinates 
to consistent system states using stochastic optimization is presented in figure
\ref{figLIFTseq}.
The main step in this
lifting procedure, which produces $u(x)$ and $a(x)$ solution profiles
with specified diffusion map coordinates, is the minimization of a
quadratic objective function given by
\begin{equation}
Obj(\Phi_2^{(i)})=\lambda_{OBJ,2}(\Phi_2^{(i)}-\Phi_2^{targ})^2
\label{objfunction}
\end{equation}
where $\lambda_{OBJ,2}$ is a weighting parameter
that controls the shape of the objective away from its minimum at
$\Phi_2^{(i)}=\Phi_2^{targ}$.
The implicit dependence of the eigenvector components
$\Phi_j^{(i)}$ (associated with a data vector ${\bf x}^{(i)}$) on ${\bf x}^{(i)}$ and all
other data vectors used in the construction of the neighbourhood matrix
${\bf K}$ makes this optimization problem challenging.
We used here, for simplicity, the method of Simulated
Annealing (SA) \cite{kirgel83,preteu92} to minimize the
objective function defined in~(\ref{objfunction}), and identify a
data vector ${\bf x}^{(targ)}$ with the target diffusion map coordinate
$\Phi_2^{targ}$. Table \ref{tab1} compares drift and
diffusion coefficient values (at $\Phi_2=-0.5$) estimated by multiple
short simulation bursts, initialized
using this lifting protocol, with values computed from a long time simulation.
The agreement is very good.
%
%
%
%

%
%
%

\begin{table}
\caption{Drift and diffusion coefficients at diffusion map coordinate
 value $\Phi_2^*=-0.5$  estimated by finding occurrences of $\Phi_2^*$ in a long
  simulation (database) and by direct initialization at this value (lifting).}
\centerline{
\begin{tabular}{|c|c|c|} \hline
 Method & $\mu(\Phi_2=\Phi_2^*)$  & $D(\Phi_2=\Phi_2^*)$ \\
\hline
Database & $-3.51\times10^{-3}$ &  $3.30\times10^{-3}$ \\
\hline
Lifting  & $-3.39\times10^{-3}$ &  $3.20\times10^{-3}$ \\
\hline
\end{tabular}}
\label{tab1}
\end{table}

Given these protocols, the entire effective potential and, more generally, all the
bifurcation/switching time computations performed above with our empirical $V$ coarse
variable can be repeated with the $\Phi_2$ variable.
``Coarse-grained smoothness" in the diffusion map coordinates can be exploited to guide
the efficient exploration of the effective potential surface \cite{frehum06}.


\section{Discussion and Conclusion}
In this paper we have revisited a stochastic spatio-temporal pattern forming system originally used
to model activity in the cortex.
Previous work~\cite{lailon01} investigated the effects of changing
two parameters (spike frequency adaptation strength and noise intensity) on the dynamics of a 
spatially-localized
``bump'' of neural activity.
Here, we have concentrated on deriving and using an effective low-dimensional description of the
system in terms of a single scalar coarse variable.
This variable was assumed to satisfy an (unknown) effective Langevin equation, so that
its probability density satisfies an effective Fokker-Planck equation.
By appropriate processing of the results of a long simulation, or --- more efficiently --- by
appropriately initializing short bursts of simulation of the full system,
we estimated the drift and diffusion functions that appear in the
Fokker-Planck equation, and thus determined an ``effective potential'' for the system on demand.
In the second part of the paper we showed how a similar analysis can be performed using a variable
that is extracted by data-mining a sufficiently long system simulation in an automated fashion,
using the diffusion map approach~\cite{nadlaf05,nadlaf06}.
Such an approach could prove particularly useful for other systems, 
since in general a good coordinate (or coordinates)
for a low-dimensional description of a complex system may not be known, nor easy to guess.
The results obtained using this new variable
showed good correspondence with those obtained in the first part of the paper.

Considering more general systems,
it is known that in certain limiting cases (e.g. weak coupling) it is possible to explicitly
obtain accurate reduced descriptions (e.g. phase oscillator equations) for large coupled
neural systems~\cite{ermpas01,izh00}.
Clearly, if such equations can be derived, they are much easier to use than the approach
presented here; our approach is intended for cases where we believe the reduction is possible,
but cannot be explicitly performed.
One of the most challenging tests for a reduction --- as we move away from the conditions where
we can guarantee its validity analytically --- is to test, on line, whether it is accurate, and
--- importantly --- whether a {\em different level} reduction, with more (or even possibly with fewer)
variables is in order.
In our case, this would correspond to devising tests to suggest, on
line, that {\em more than one}
coarse-grained variable must take part in our effective Langevin equation. The slight ``thickness"
of the line in figure~\ref{fig:coordMAP} could be an indication that, for the parameter values used, more than one variable may be necessary
in our model.
Devising such tests --- in effect, testing the hypothesis that the data are locally consistent
with a particular model, e.g. a scalar Langevin equation --- is the subject of ongoing research across
several disciplines (statistics, financial mathematics). Integrating such tools in a multiscale
simulation framework is only just starting.

Even though the model discussed here is rather abstract (from a computational neuroscience point of view) it adds to the
body of work demonstrating that it is possible to obtain useful low-dimensional descriptions of complex systems,
both in neural modelling~\cite{lai06} and elsewhere~\cite{erbkev06,kevgea03,srikev05}.
These descriptions not only
provide insight into the fundamental dynamics, but enable one to simulate and analyse such systems in an
efficient manner.

Coarse-graining large scale, faithful neural network computations constitutes an important subject of
intense current research (see, for example, the probability density approach~\cite{caitao04,kni00,omukni00}).
Multiscale, coarse-graining numerical algorithms such as the one demonstrated here
have an important role to play in elucidating the types of behavior possible for such
models and their parametric dependence.

\ack We are pleased to acknowledge discussions with Professor R. Coifman and
members of his group on data mining and diffusion maps. 

\appendix
\section*{Appendix}
\setcounter{section}{1}


We run a simulation for 30,000 time units and sample the $u(x),a(x)$ profile every 8 time units,
giving $N=3,\!750$ data points. We have discretised space with $M=100$ points, so every data point is actually
a vector of 200 values.
The components of the simulation data vector ${\bf x}$ for our problem
are determined by the $2M$ nodal values of the dependent variables $u_i$ and $a_i$
(appearing in (\ref{eq:dudtA})-(\ref{eq:dadtA})) at a particular time step. We normalize each
$u_i$ profile by its maximum value (at node $i_{max}$) and normalize the
$a_i$ profile by its nodal value at the corresponding location
($a_{i_{max}}$). Finally, we shift both solution profiles by a fixed
number of nodes such that the alignment between $u_i$ profiles is
maximized (in a least squares sense) thereby factoring out
translations.  We define a pairwise similarity (``neighbourhood'') matrix ${\bf K}$
between a representative sample of these data vectors (collected
over the course of a simulation run) as follows
\begin{equation}
 K_{i,j}=K\left({\bf x}_i,{\bf x}_j\right)=\exp\left[-\left(\frac{\norm{{\bf x}_i-{\bf x}_j}}{\sigma}\right)^2\right],
\label{kernel}
\end{equation}
where $\sigma$ is a parameter that defines the size
of the local neighbourhood surrounding each high-dimensional point. Defining the diagonal normalization matrix $D_{i,i}=\sum_j K_{i,j}$,
we construct the Markovian matrix ${\bf M=D}^{-1}{\bf K}$. A few top
eigenvectors of the matrix ${\bf M}$ are used here as a low dimensional
representation of the simulation data. Components in the eigenvector provide
the low dimensional coordinate for each simulation data
vector. The diffusion map distance is equal to Euclidean distance in
{\em the
diffusion map space}. In many applications the spectrum of the matrix
${\bf M}$ possesses a spectral gap, and this diffusion distance may be
approximated by using only a few top eigenvectors. This 
approximation has been shown to be optimal under a particular mean
squared error criterion \cite{nadlaf05}. 

The matrix ${\bf M}$ is adjoint to the symmetric matrix ${\bf M_s}$
defined as follows
\begin{equation}
 {\bf M_s}={\bf D}^{1/2}{\bf M D}^{-1/2}={\bf D}^{-1/2}{\bf K D}^{-1/2}.
\label{adjointM}
\end{equation}
${\bf M}$ and ${\bf M_s}$ share the same eigenvalues; eigenvectors
${\bf\Phi}_j$ of ${\bf M}$ are related to those of ${\bf M}_s$, denoted ${\bf\Psi}_j$, as follows
\begin{equation}
 {\bf\Phi}_j={\bf D}^{-1/2}{\bf\Psi}_j
\label{eigenvecM+Ms}
\end{equation}
We compute the largest eigenvalues and eigenvectors of the symmetric
matrix ${\bf M_s}$ and use~(\ref{eigenvecM+Ms}) to evaluate the
eigenvectors of the Markov matrix ${\bf M}$. We note that the top eigenvector
${\bf\Psi}_1$ of ${\bf M_s}$, corresponding to the largest eigenvalue
$\lambda_1=1$, has components equal to the diagonal entries of ${\bf
  D}^{1/2}$; it follows from~(\ref{eigenvecM+Ms}) that the top,
{\em trivial}, eigenvector ${\bf\Phi}_1$ of ${\bf M}$ consists entirely
of ones.

\end{document}